\documentclass[a4paper,12pt, reqno]{amsart}
\usepackage{latexsym}
\usepackage{amssymb,amsfonts,amsmath,mathrsfs}
\begin{document}
\title{On the mean values of $L$-functions in orthogonal and symplectic families}
\author{H. M. BUI, J. P. KEATING}
\address{School of Mathematics, University of Bristol, Bristol, BS8 1TW}
\email{hm.bui@bristol.ac.uk, j.p.keating@bristol.ac.uk}

\begin{abstract}
Hybrid Euler-Hadamard products have previously been studied for the Riemann zeta function on its critical line and for Dirichlet $L$-functions, in the context of the calculation of moments and connections with Random Matrix Theory. According to the Katz-Sarnak classification, these are believed to represent families of $L$-function with unitary symmetry. We here extend the formalism to families with orthogonal \& symplectic symmetry. Specifically, we establish formulae for real quadratic Dirichlet $L$-functions and for the $L$-functions associated with primitive Hecke eigenforms of weight 2 in terms of partial Euler and Hadamard products. We then prove asymptotic formulae for some moments of these partial products and make general conjectures based on results for the moments of characteristic polynomials of random matrices.
\end{abstract}
\maketitle

\section{Introduction}

Estimating the moments of families of $L$-functions is a central theme in number theory. These have extensive applications, such as to bounding the order of $L$-functions and to proving results relating to non-vanishing. Moreover, the problems are interesting in their own right, since, according to Katz and Sarnak's philosophy [\textbf{\ref{KS1}}], they are expected to illustrate the symmetry of the families.  

The most well-understood mean values are the moments of the Riemann zeta function. It has long been conjectured that for every $k\geq0$, as $T\rightarrow\infty$, there is a constant $c_k$ such that
\begin{equation}\label{1}
I_{k}(T):=\frac{1}{T}\int_{0}^{T}|\zeta({\scriptstyle{\frac{1}{2}}}+it)|^{2k}dt\sim c_{k}(\log T)^{k^2}.
\end{equation}
The second moment was established by Hardy and Littlewood [\textbf{\ref{HL}}],
\begin{displaymath}
I_{1}(T)\sim\log T,
\end{displaymath} 
and Ingham [\textbf{\ref{I}}] gave an asymptotic formula for the fourth moment,
\begin{displaymath}
I_{2}(T)\sim\frac{1}{2\pi^2}(\log T)^4.
\end{displaymath}
No other mean values of the zeta function have been proved. Conrey and Ghosh [\textbf{\ref{CG1}}] cast \eqref{1} in a more precise form, namely, there should be a factorization
\begin{displaymath}
c_{k}=\frac{a_{k,U}g_{k,U}}{\Gamma(k^2+1)},
\end{displaymath}
where
\begin{equation}
a_{k,U}=\prod_{p}\bigg[\bigg(1-\frac{1}{p}\bigg)^{k^2}\sum_{m\geq0}\frac{d_{k}(p^m)^2}{p^m}\bigg],
\end{equation}
and $g_{k,U}$ is an integer when $k$ is an integer. The classical results of Hardy-Littlewood and Ingham imply $g_{1,U}=1$ and $g_{2,U}=2$. Conrey and Ghosh [\textbf{\ref{CG2}}], and then later Conrey and Gonek [\textbf{\ref{CGo}}] used Dirichlet polynomial techniques to conjecture that $g_{3,U}=42$ and $g_{4,U}=24024$, respectively. Their method also reproduced the previous values of $g_{k,U}$ but does not give conjectural values for $g_{k,U}$ for $k>4$. However, we do have a lower bound, and, assuming the Riemann Hypothesis, an upper bound for $I_{k}(T)$, for every positive real number $k$, that are consistent with \eqref{1}. See [\textbf{\ref{Sound2}}] for a brief survey and significant new results in this direction.

For the family of Dirichlet $L$-functions $L(s,\chi_d)$ with real primitive Dirichlet characters $\chi_d$ modulo $d$, Jutila [\textbf{\ref{Jutila}}] proved that as $D\rightarrow\infty$,
\begin{equation}
\frac{1}{D}\sum_{0<d\leq D}L({\scriptstyle{\frac{1}{2}}},\chi_d)\sim\frac{1}{4\zeta(2)}\prod_{p}\bigg(1-\frac{1}{p(p+1)}\bigg)\log D,
\end{equation}
and
\begin{equation}
\frac{1}{D}\sum_{0<d\leq D}L({\scriptstyle{\frac{1}{2}}},\chi_d)^{2}\sim\frac{1}{48\zeta(2)}\prod_{p}\bigg(1-\frac{4p^2-3p+1}{p^3(p+1)}\bigg)(\log D)^3.
\end{equation}
If we restrict $d$ to be odd, positive and square-free integers so that $\chi_{8d}$ are real, primitive characters with conductors $8d$ and with $\chi_{8d}(-1)=1$, then the corresponding results are
\begin{equation}\label{2}
\frac{1}{D^{*}}\sum_{0<d\leq D}{\!\!\!\!}^{\flat}\ L({\scriptstyle{\frac{1}{2}}},\chi_{8d})\sim\frac{1}{4}\prod_{p\geq3}\bigg(1-\frac{1}{p(p+1)}\bigg)\log D,
\end{equation}
and
\begin{equation}\label{3}
\frac{1}{D^{*}}\sum_{0<d\leq D}{\!\!\!\!}^{\flat}\ L({\scriptstyle{\frac{1}{2}}},\chi_{8d})^{2}\sim\frac{1}{192}\prod_{p\geq3}\bigg(1-\frac{4p^2-3p+1}{p^3(p+1)}\bigg)(\log D)^3,
\end{equation}
where the sum $\sum^{\flat}$ over $d$ indicates that $d$ is odd and square-free, and $D^{*}$ is the number of such $d$ in $(0,D]$.

Recently, Soundararajan [\textbf{\ref{Sound1}}] showed that
\begin{eqnarray}\label{4}
\frac{1}{D^{*}}\sum_{0<d\leq D}{\!\!\!\!}^{\flat}\ L({\scriptstyle{\frac{1}{2}}},\chi_{8d})^{3}&\sim&\nonumber\\
&&\!\!\!\!\!\!\!\!\!\!\!\!\!\!\!\!\!\!\!\!\!\!\!\!\!\!\!\!\!\!\!\!\!\!\!\!\!\!\!\!\!\!\!\!\!\!\!\!\!\!\!\!\frac{1}{184320}\prod_{p\geq3}\bigg(1-\frac{12p^5-23p^4+23p^3-15p^2+6p-1}{p^6(p+1)}\bigg)(\log D)^6.
\end{eqnarray}
Soundararajan also gave a conjecture for the fourth moment, but the higher moments are not well-understood.

Another family of $L$-functions which has been considered by a number of authors arises from $S_{2}^{*}(q)$, the set of primitive Hecke eigenforms of weight $2$ and level $q$ ($q$ prime). For $f(z)\in S_{2}^{*}(q)$, $f$ has a Fourier expansion
\begin{displaymath}
f(z)=\sum_{n\geq1}n^{1/2}\lambda_{f}(n)e(nz),
\end{displaymath}
where the normalization is such that $\lambda_{f}(1)=1$. The $L$-function associated to $f$ has an Euler product
\begin{displaymath}
L(f,s)=\sum_{n\geq1}\frac{\lambda_{f}(n)}{n^s}=\bigg(1-\frac{\lambda_{f}(q)}{q^s}\bigg)^{-1}\prod_{p\ne q}\bigg(1-\frac{\lambda_{f}(p)}{p^s}+\frac{1}{p^{2s}}\bigg)^{-1},
\end{displaymath} 
which we can express as
\begin{displaymath}
L(f,s)=\prod_{p}\bigg(1-\frac{\alpha_{f}(p)}{p^s}\bigg)^{-1}\bigg(1-\frac{\alpha'_{f}(p)}{p^s}\bigg)^{-1},
\end{displaymath}
where $\alpha'_{f}(q)=0$ and $\alpha'_{f}(p)=\overline{\alpha_{f}(p)}$ if $p\neq q$. The series is absolutely convergent when $\Re s>1$ and admits an analytic continuation to all of $\mathbb{C}$. The functional equation for the $L$-function is
\begin{displaymath}
\Lambda(f,s):=\bigg(\frac{\sqrt{q}}{2\pi}\bigg)^{s}\Gamma(s+{\scriptstyle{\frac{1}{2}}})L(f,s)=\epsilon_{f}\Lambda(f,1-s),
\end{displaymath}
where $\epsilon_f=q^{1/2}\lambda_{f}(q)=\pm1$. We define the harmonic average as
\begin{displaymath}
\sum_{f}{\!}^{h}\ A_{f}:=\sum_{f\in S_{2}^{*}(q)}\frac{A_{f}}{4\pi(f,f)},
\end{displaymath}
where $(f,g)$ is the Petersson inner product on the space $\Gamma_{0}(q)\backslash\mathbb{H}$.

It has been shown by Duke [\textbf{\ref{D}}], Duke, Friedlander, and Iwaniec [\textbf{\ref{DFI}}], and Iwaniec and Sarnak [\textbf{\ref{IS}}] that as $q\rightarrow\infty$,
\begin{equation}\label{13}
\sum_{f}{\!}^{h}\ L(f,{\scriptstyle{\frac{1}{2}}})\sim1,
\end{equation}
and
\begin{equation}\label{14}
\sum_{f}{\!}^{h}\ L(f,{\scriptstyle{\frac{1}{2}}})^{2}\sim\log q,
\end{equation}
and by Kowalski, Michel, and VanderKam [\textbf{\ref{KMV1}}] [\textbf{\ref{KMV2}}] that
\begin{equation}\label{15}
\sum_{f}{\!}^{h}\ L(f,{\scriptstyle{\frac{1}{2}}})^{3}\sim\frac{1}{6}(\log q)^3,
\end{equation}
and
\begin{equation}\label{16}
\sum_{f}{\!}^{h}\ L(f,{\scriptstyle{\frac{1}{2}}})^{4}\sim\frac{1}{360\zeta(2)}(\log q)^6.
\end{equation}

Following an idea of Katz and Sarnak [\textbf{\ref{KS1}}] which associates a family of $L$-functions to a corresponding symmetry group and asserts that the symmetry group governs many properties of the distribution of zeros of the $L$-functions, Conrey and Farmer [\textbf{\ref{CF}}] gave evidence that the symmetry type of a family of $L$-functions also determines the mean values of the $L$-functions at the critical point. Precisely, they conjectured that in general
\begin{equation}\label{6}
\frac{1}{\mathcal{Q}^{*}}\sum_{\substack{f\in\mathcal{F}\\c(f)\leq \mathcal{Q}}}V(L_{f}({\scriptstyle{\frac{1}{2}}}))^{k}\sim \frac{a_{k}g_{k}}{\Gamma(1+B(k))}(\log\mathcal{Q})^{B(k)},
\end{equation}
for some $a_k$, $g_k$, and $B(k)$, where the $L$-functions are normalized to have a functional equation $s\leftrightarrow1-s$, so $\frac{1}{2}$ is the central point; the family $\mathcal{F}$ is considered to be partially ordered by the conductor $c(f)$, and $\mathcal{Q}^{*}$ is the number of elements with $c(f)\leq\mathcal{Q}$; and $V(z)$ is chosen depending on the symmetry type of the family. Conrey and Farmer observed that $g_k$ and $B(k)$ depend only on the symmetry type of the family, while $a_k$ depends on the family itself and is computable in any specific case. In this shape, $V(z)=|z|^2$ for unitary symmetry and $V(z)=z$ for symplectic or orthogonal symmetry. The Riemann zeta function is viewed as forming its own unitary family and has $B(k)=k^2$. The family of Dirichlet $L$-functions $L(s,\chi_{8d})$, where $d$ are odd, positive and square-free integers is believed to have symplectic symmetry. In this case, $B(k)=k(k+1)/2$ and
\begin{equation}
a_{k,Sp}=2^{-k(k+2)/2}\prod_{p\geq3}\frac{(1-\frac{1}{p})^{k(k+1)/2}}{1+\frac{1}{p}}\bigg(\frac{(1-\frac{1}{\sqrt{p}})^{-k}+(1+\frac{1}{\sqrt{p}})^{-k}}{2}+\frac{1}{p}\bigg).
\end{equation}
\textit{Remark}. An equivalent form of $a_{k,Sp}$ is
\begin{equation}
a_{k,Sp}=2^{-k(k+2)/2}\prod_{p\geq3}\bigg[\bigg(1-\frac{1}{p}\bigg)^{k(k+1)/2}\bigg(1+\sum_{m=1}^{\infty}\frac{d_{k}(p^{2m})}{p^{m-1}(p+1)}\bigg)\bigg].
\end{equation}
The family of $L$-functions $L(f,s)$, where $f\in S_{2}^{*}(q)$, is conjectured to be included in the orthogonal symmetry type with $B(k)=k(k-1)/2$. The constant $a_{k,O}$ associated to this family can also be determined explicitly,
\begin{equation}
a_{k,O}=2^{k/2}\zeta(2)^{k}\prod_{p}\bigg[\bigg(1-\frac{1}{p}\bigg)^{k(k+1)/2}\sum_{m\geq 0}\frac{d_{k}(p^{2m})}{p^{m}(1+1/p)^{2m}}\bigg(\binom{2m}{m}-\binom{2m}{m-1}\bigg)\bigg].
\end{equation}

While the conjecture \eqref{6} is verified only for some small values of $k$, all the parameters in the formula are fairly well understood, except for the constant $g_k$. This is the motivation for the random matrix model introduced by Keating and Snaith [\textbf{\ref{KS2}}], in which statistical properties of the Riemann zeta function are related to those of the
characteristic polynomials of large random matrices. Specifically, let $U$ be an $N\times N$ unitary matrix. If we denote the eigenvalues of $U$ by $e^{i\theta_n}$, the characteristic polynomial of $U$ is
\begin{displaymath}
Z_{N}(U,\theta)=\prod_{n=1}^{N}(1-e^{i(\theta_{n}-\theta)}).
\end{displaymath}
Keating and Snaith then proved that as $N\rightarrow\infty$,
\begin{equation}\label{7}
\mathbb{E}_{N}\big[|Z_{N}(U,\theta)|^{2k}\big]\sim
\frac{G^{2}(k+1)}{G(2k+1)}N^{k^2},
\end{equation}
where the expectation value is computed with respect to Haar measure on $U(N)$, and $G(z)$ is Barnes' $G$-function. Equating the mean density of the eigenphases $\theta_n$ to the mean density of the zeros of the Riemann zeta function corresponds to the identification $N\sim\log T$, and hence $N^{k^2}$ gives the right order for the 2$k^{\textrm{th}}$ moment of the zeta function. Also, based on the fact that $G^{2}(k+1)/G(2k+1)=g_{k,U}/\Gamma(k^2+1)$ for $k=1,2$, and, conjecturally, for $k=3,4$, Keating and Snaith conjectured that this formula holds in general, that is
\newtheorem{conj}{Conjecture}\begin{conj}
For $k$ fixed with $\Re k>-1/2$, as $T\rightarrow \infty$, 
\begin{displaymath}
\frac{1}{T}\int_{0}^{T}|\zeta({\scriptstyle{\frac{1}{2}}}+it)|^{2k}dt\sim a_{k,U}\frac{G^{2}(k+1)}{G(2k+1)}(\log T)^{k^2}.
\end{displaymath}
\end{conj}
Extending their results to the symplectic and orthogonal groups, Keating and Snaith [\textbf{\ref{KS3}}] then derived conjectures for the mean values of $L$-functions in symplectic and orthogonal families.
\begin{conj}
For $k$ fixed with $\Re k\geq0$, as $D\rightarrow \infty$,
\begin{displaymath}
\frac{1}{D^{*}}\sum_{0<d\leq D}{\!\!\!\!}^{\flat}\ L({\scriptstyle{\frac{1}{2}}},\chi_{8d})^{k}\sim a_{k,Sp}\frac{G(k+1)\sqrt{\Gamma(k+1)}}{\sqrt{G(2k+1)\Gamma(2k+1)}}(\log D)^{k(k+1)/2}.
\end{displaymath}
\end{conj}
\begin{conj}
For $k$ fixed with $\Re k\geq0$, as $q\rightarrow \infty$,
\begin{displaymath}
\sum_{f}{\!}^{h}\ L(f,{\scriptstyle{\frac{1}{2}}})^{k}\sim a_{k,O}\frac{G(k+1)\sqrt{\Gamma(2k+1)}}{2\sqrt{G(2k+1)\Gamma(k+1)}}(\log q)^{k(k-1)/2}.
\end{displaymath}
\end{conj}
These conjectures are discussed in detail and extended to include all lower order terms in the asymptotic series representing the moments in [\textbf{\ref{CFKRS}}].

However, the drawback of the random matrix model is the absence of the arithmetical factors $a_k$. This can be obtained from number-theoretical considerations [\textbf{\ref{CF}}], [\textbf{\ref{CFKRS}}], but then the random-matrix contribution appears mysteriously. The question is how to treat the arithmetical and random-matrix aspects on an equal footing.

Recently, Gonek, Hughes and Keating [\textbf{\ref{GHK}}] proved, using a smoothed form of the explicit formula of Bombieri and Hejhal [\textbf{\ref{BH}}], that one can approximate the Riemann zeta function at a height $t$ on the critical line as a partial Euler product multiplied by a partial Hadamard product over the nontrivial zeros close to $1/2 +it$. This suggests a statistical model for the zeta function in which the primes are incorporated in a natural way. The
value distribution of the product over zeros is expected to be modelled by the characteristic polynomial of a large random unitary matrix, because it involves only local information about the zeros. Conjecturing the moments of this product using random matrix theory, calculating the moments of the product over the primes rigorously and making an assumption (which can be proved in certain particular cases) about the independence of the two products, Gonek, Hughes and Keating then reproduced the conjecture for the $2k^{\textrm{th}}$ moment of the zeta function first put forward by Keating and Snaith in [\textbf{\ref{KS2}}].

In our previous paper [\textbf{\ref{BK}}], we extended this approach to the $2k^{\textrm{th}}$ power moment of Dirichlet
\textit{L}-functions $L(s,\chi)$ at the centre of the critical strip ($s=1/2$), where the average is over all primitive characters $\chi$ (mod $q$). These $L$-functions form a unitary family, so the results are similar to the zeta function case. Here we show that the model introduced in [\textbf{\ref{GHK}}] can be adapted to $L$-functions with symplectic or orthogonal symmetry.

\section{Statement of results}

\subsection{Symplectic family.} 

The Euler-Hadamard product for the quadratic Dirichlet $L$-functions is stated in [\textbf{\ref{BK}}] (cf. Theorem 1 there).\vspace{0.5cm}\\
\textbf{Theorem S1} \textit{Let $u(x)$ be a real, non-negative, $C^{\infty}$ function with mass $1$ and
compact support on $[e^{1-1/X},e]$. Set
\begin{displaymath}
U(z)=\int_{0}^{\infty}u(x)E_{1}(z\log x)dx,
\end{displaymath}
where $E_{1}(z)$ is the exponential integral $\int_{z}^{\infty}e^{-w}/wdw$. Let $X\geq2$ be a real parameter.
Then for $d$ positive, odd and square-free,
\begin{equation}
L({\scriptstyle{\frac{1}{2}}},\chi_{8d})=P_{X}(\chi_{8d})Z_{X}(\chi_{8d})(1+O((\log X)^{-2})),
\end{equation}
where
\begin{equation}
P_{X}(\chi_{8d})=\exp\bigg( \sum_{n\leq X}\frac{\Lambda(n)\chi_{8d}(n)}{n^{1/2}\log n}\bigg),
\end{equation}
$\Lambda(n)$ is von Mangoldt's function, and
\begin{equation}
Z_{X}(\chi_{8d})=\exp\bigg(-\sum_{\rho}U(({\scriptstyle{\frac{1}{2}}}-\rho)\log X)\bigg),
\end{equation}
where the sum is over the nontrivial zeros $\rho$ of $L(s,\chi_{8d})$. The constant implied by the $O$-term is absolute.}\\

We calculate the moments of $P_{X}(\chi_{8d})$ rigorously and establish:\vspace{0.5cm}\\
\textbf{Theorem S2} \textit{Let $\delta>0$ and $k\geq0$ fixed. Suppose that $X,D\rightarrow\infty$ with $X\ll(\log D)^{2-\delta}$, then
\begin{displaymath}
\frac{1}{D^{*}}\sum_{0<d\leq D}{\!\!\!\!}^{\flat}\ P_{X}(\chi_{8d})^{k}\sim a_{k,Sp}(e^{\gamma}\log X)^{k(k+1)/2}.
\end{displaymath}}

We note from Theorem S1 that $L(1/2,\chi_{8d})P_{X}(\chi_{8d})^{-1}=Z_{X}(\chi_{8d})(1+o(1))$. This allows us to derive the first two moments of $Z_X$.\vspace{0.5cm}\\
\textbf{Theorem S3} \textit{Let $\delta>0$. For $X,D\rightarrow\infty$ with $X\ll(\log D)^{2-\delta}$, we have
\begin{displaymath}
\frac{1}{D^{*}}\sum_{0<d\leq D}{\!\!\!\!}^{\flat}\ L({\scriptstyle{\frac{1}{2}}},\chi_{8d})P_{X}(\chi_{8d})^{-1}\sim\frac{1}{\sqrt{2}}\frac{\log D}{e^{\gamma}\log X}.
\end{displaymath}}\vspace{0.05cm}\\
\textbf{Theorem S4} \textit{Let $\delta>0$. For $X,D\rightarrow\infty$ with $X\ll(\log D)^{2-\delta}$, we have
\begin{displaymath}
\frac{1}{D^{*}}\sum_{0<d\leq D}{\!\!\!\!}^{\flat}\ L({\scriptstyle{\frac{1}{2}}},\chi_{8d})^{2}P_{X}(\chi_{8d})^{-2}\sim\frac{1}{12}\bigg(\frac{\log D}{e^{\gamma}\log X}\bigg)^3.
\end{displaymath}}

We remark that these asymptotic results coincide precisely with the corresponding formulae for the moments of the characteristic polynomials of random symplectic matrices if we identify $\log D/e^\gamma\log X$ with the matrix size. This is expected, on the basis of the heuristic arguments given by Gonek, Hughes and Keating [\textbf{\ref{GHK}}]. We hence conjecture that this holds in general. \vspace{0.5cm}\\
\textbf{Conjecture S1} \textit{Let $k\geq0$ be any real number. Suppose that $X$ and $D\rightarrow\infty$ with $X\ll(\log D)^{2-\epsilon}$, then
\begin{displaymath}
\frac{1}{D^{*}}\sum_{0<d\leq D}{\!\!\!\!}^{\flat}\ Z_{X}(\chi_{8d})^{k}\sim\frac{G(k+1)\sqrt{\Gamma(k+1)}}{\sqrt{G(2k+1)\Gamma(2k+1)}}\bigg(\frac{\log D}{e^\gamma\log X}\bigg)^{k(k+1)/2}.
\end{displaymath}}

Combining the formulae for the first ($k=1$) and second ($k=2$) moments, \eqref{2} and \eqref{3}, with Theorem S2, Theorem S3 and Theorem S4, we see that, at least for the cases $k=1$ and $k=2$, when $X$ is not too large relative to $D$, the $k^{\textrm{th}}$ moment of $L(1/2,\chi_{8d})$ is asymptotic to the product of the moments of $P_{X}(\chi_{8d})$ and $Z_{X}(\chi_{8d})$. We believe that this is true in general:\vspace{0.5cm}\\
\textbf{Conjecture S2} \textit{Let $k\geq0$ be any real number. Suppose that $X$ and $D\rightarrow\infty$ with $X\ll(\log D)^{2-\epsilon}$, then
\begin{displaymath}
\frac{1}{D^{*}}\sum_{0<d\leq D}{\!\!\!\!}^{\flat}\ L({\scriptstyle{\frac{1}{2}}},\chi_{8d})^{k}\sim\bigg(\frac{1}{D^{*}}\sum_{0<d\leq D}{\!\!\!\!}^{\flat}\ P_{X}(\chi_{8d})^{k}\bigg)\bigg(\frac{1}{D^{*}}\sum_{0<d\leq D}{\!\!\!\!}^{\flat}\ Z_{X}(\chi_{8d})^{k}\bigg).
\end{displaymath}}

\subsection{Orthogonal family.} 

The Euler-Hadamard product for the $L$-functions associated with $f\in S_{2}^{*}(q)$ can be proved by following step-by-step the corresponding calculation in [\textbf{\ref{BK}}] (cf. Theorem 1 there).\vspace{0.5cm}\\
\textbf{Theorem O1} \textit{Let $u(x)$ be a real, non-negative, $C^{\infty}$ function with mass $1$ and
compact support on $[e^{1-1/X},e]$. Set
\begin{displaymath}
U(z)=\int_{0}^{\infty}u(x)E_{1}(z\log x)dx.
\end{displaymath}
Then for $f\in S_{2}^{*}(q)$,
\begin{equation}
L(f,{\scriptstyle{\frac{1}{2}}})=P_{X}(f)Z_{X}(f)(1+O((\log X)^{-2})),
\end{equation}
where
\begin{equation}
P_{X}(f)=\exp\bigg(\sum_{n\leq X}\frac{\Lambda_{f}(n)}{n^{1/2}\log n}\bigg),
\end{equation}
and
\begin{equation}
Z_{X}(f)=\exp\bigg(-\sum_{\rho}U(({\scriptstyle{\frac{1}{2}}}-\rho)\log X)\bigg),
\end{equation}
where
\begin{displaymath}
\Lambda_{f}(n)=\left\{ \begin{array}{ll}
(\alpha_{f}(p)^{j}+\alpha'_{f}(p)^{j})\log p &\qquad \textrm{if $n=p^j$}\\
0 & \qquad\textrm{otherwise,} 
\end{array} \right.
\end{displaymath}
and the sum in $Z_X$ is over the nontrivial zeros $\rho$ of $L(f,s)$. The constant implied by the $O$-term is absolute.}\\

The mean values of $P_X(f)$ can be established rigorously.\vspace{0.5cm}\\
\textbf{Theorem O2} \textit{Let $\delta>0$ and $k\geq0$ fixed. Suppose that $X,q\rightarrow\infty$ with $X\ll(\log q)^{2-\delta}$, then
\begin{displaymath}
\sum_{f}{\!}^{h}\ P_{X}(f)^{k}\sim a_{k,O}(e^{\gamma}\log X)^{k(k-1)/2}.
\end{displaymath}}

We note from Theorem O1 that $L(f,1/2)P_{X}(f)^{-1}=Z_{X}(f)(1+o(1))$. This allows us to derive the first four harmonic moments of $Z_X$.\vspace{0.5cm}\\
\textbf{Theorem O3} \textit{Let $\delta>0$. For $X,q\rightarrow\infty$ with $X\ll(\log q)^{2-\delta}$, we have
\begin{displaymath}
\sum_{f}{\!}^{h}\ L(f,{\scriptstyle{\frac{1}{2}}})P_{X}(f)^{-1}\sim\frac{1}{\sqrt{2}}.
\end{displaymath}}\vspace{0.05cm}\\
\textbf{Theorem O4} \textit{Let $\delta>0$. For $X,q\rightarrow\infty$ with $X\ll(\log q)^{2-\delta}$, we have
\begin{displaymath}
\sum_{f}{\!}^{h}\ L(f,{\scriptstyle{\frac{1}{2}}})^{2}P_{X}(f)^{-2}\sim \frac{1}{2}\frac{\log q}{e^{\gamma}\log X}.
\end{displaymath}}\vspace{0.05cm}\\
\textbf{Theorem O5} \textit{Let $\delta>0$. For $X,q\rightarrow\infty$ with $X\ll(\log q)^{2-\delta}$, we have
\begin{displaymath}
\sum_{f}{\!}^{h}\ L(f,{\scriptstyle{\frac{1}{2}}})^{3}P_{X}(f)^{-3}\sim \frac{1}{12\sqrt{2}}\bigg(\frac{\log q}{e^{\gamma}\log X}\bigg)^3.
\end{displaymath}}\vspace{0.05cm}\\
\textbf{Theorem O6} \textit{Let $\delta>0$. For $X,q\rightarrow\infty$ with $X\ll(\log q)^{2-\delta}$, we have
\begin{displaymath}
\sum_{f}{\!}^{h}\ L(f,{\scriptstyle{\frac{1}{2}}})^{4}P_{X}(f)^{-4}\sim \frac{1}{1440}\bigg(\frac{\log q}{e^{\gamma}\log X}\bigg)^6.
\end{displaymath}}

As in the symplectic case, these formulae coincide with the moments of the characteristic polynomials of random orthogonal matrices if $\log q/e^\gamma\log X$ is identified with the matrix size. We expect this to hold in general.\vspace{0.5cm}\\
\textbf{Conjecture O1} \textit{Let $k\geq0$ be any real number. Suppose that $X$ and $q\rightarrow\infty$ with $X\ll(\log q)^{2-\epsilon}$, then
\begin{displaymath}
\sum_{f}{\!}^{h}\ Z_{X}(f)^{k}\sim\frac{G(k+1)\sqrt{\Gamma(2k+1)}}{2\sqrt{G(2k+1)\Gamma(k+1)}}\bigg(\frac{\log q}{e^\gamma\log X}\bigg)^{k(k-1)/2}.
\end{displaymath}}

Combining the results \eqref{13}, \eqref{14}, \eqref{15}, and \eqref{16}, with Theorem O2-O6, we see that, at least for the cases $k=1,2,3$ and $4$, when $X$ is not too large relative to $q$, the harmonic $k^{\textrm{th}}$ moment of $L(f,1/2)$ is asymptotic to the product of the harmonic moments of $P_{X}(f)$ and $Z_{X}(f)$. We believe that this is true in general:\vspace{0.5cm}\\
\textbf{Conjecture O2} \textit{Let $k\geq0$ be any real number. Suppose that $X$ and $q\rightarrow\infty$ with $X\ll(\log q)^{2-\epsilon}$, then
\begin{displaymath}
\sum_{f}{\!}^{h}\ L(f,{\scriptstyle{\frac{1}{2}}})^{k}\sim\bigg(\sum_{f}{\!}^{h}\ P_{X}(f)^{k}\bigg)\bigg(\sum_{f}{\!}^{h}\ Z_{X}(f)^{k}\bigg).
\end{displaymath}}

The paper is organized as follows: In the next section we consider the mean values of real quadratic Dirichlet $L$-functions; The final section is devoted to the family of $L$-functions with orthogonal symmetry type.

\section{Symplectic family}

\subsection{Proof of Theorem S2.}

We record a lemma from [\textbf{\ref{BK}}] (cf. Lemma 3).
\newtheorem{lemm}{Lemma}\begin{lemm}
Let
\begin{displaymath}
P_{X}(s,\chi_{8d})=\exp\bigg( \sum_{n\leq X}\frac{\Lambda(n)\chi_{8d}(n)}{n^{s}\log n}\bigg),
\end{displaymath} 
so $P_{X}(\chi_{8d})=P_{X}(1/2,\chi_{8d})$, and let $P_{X}^{*}(\chi_{8d})=P_{X}^{*}(1/2,\chi_{8d})$, where
\begin{displaymath}
P_{X}^{*}(s,\chi_{8d})=\prod_{p\leq X}\bigg( 1-\frac{\chi_{8d}(p)}{p^s}\bigg) ^{-1}\prod_{\sqrt{X}<p\leq X}\bigg( 1+\frac{\chi_{8d}(p)^2}{2p^{2s}}\bigg) ^{-1}.
\end{displaymath} 
Then for any $k\in\mathbb{R}$ we have
\begin{displaymath}
P_{X}(s,\chi_{8d})^{k}=P_{X}^{*}(s,\chi_{8d})^{k}(1+O_{k}(X^{-1/6+\epsilon})),
\end{displaymath}
uniformly for $\sigma\geq1/2$.
\end{lemm}
\textit{Remark}. We note that the $O$-term has been made a power in $X$. This can be done easily with a little care when dealing with the error term. See the proof of Lemma 6 below.\\

The next lemma is standard.
\begin{lemm}
For $m$ odd
\begin{displaymath}
\sum_{0<d\leq D}{\!\!\!\!}^{\flat}\ \chi_{8d}(m^2)=\frac{2Da(m)}{3\zeta(2)}+O(D^{1/2}m^\epsilon),
\end{displaymath}
where $a(m)=\prod_{p|m}\big(1+\frac{1}{p}\big)^{-1}.$
\end{lemm}
\begin{proof}
We have
\begin{eqnarray*}
\sum_{0<d\leq D}{\!\!\!\!}^{\flat}\ \chi_{8d}(m^2)&=&\sum_{\substack{0<d\leq D\\(d,2m)=1}}\sum_{l^2|d}\mu(l)=\sum_{\substack{0<l\leq\sqrt{D}\\(l,2m)=1}}\mu(l)\sum_{\substack{d\leq D/l^2\\(d,2m)=1}}1\\
&=&\frac{\varphi(2m)}{2m}D\sum_{\substack{0<l\leq\sqrt{D}\\(l,2m)=1}}\frac{\mu(l)}{l^2}+O(D^{1/2}m^\epsilon)
\end{eqnarray*}
We observe that
\begin{displaymath}
\sum_{\substack{0<l\leq\sqrt{D}\\(l,2m)=1}}\frac{\mu(l)}{l^2}=\frac{1}{\zeta(2)}\prod_{p|2m}\bigg(1-\frac{1}{p^2}\bigg)^{-1}(1+O(D^{-1/2})).
\end{displaymath}
The lemma follows.
\end{proof}
\textit{Remark}. In particular, for $m=1$, we obtain $D^{*}\sim 2D/3\zeta(2)$.\\

We now proceed with the proof of the theorem. We write $P_{X}^{*}(s,\chi_{8d})^{k}$ as a Dirichlet series
\begin{equation}
\sum_{n=1}^{\infty}\frac{\alpha_{k}(n)\chi_{8d}(n)}{n^s}=\prod_{p\leq X}\bigg( 1-\frac{\chi_{8d}(p)}{p^s}\bigg) ^{-k}\prod_{\sqrt{X}<p\leq X}\bigg( 1+\frac{\chi_{8d}(p)^2}{2p^{2s}}\bigg) ^{-k}.
\end{equation}
We note that $\alpha_{k}(n)\in\mathbb{R}$, and if we denote by $S(X)$ the set of $X$-smooth numbers, i.e.
\begin{displaymath}
S(X)=\{n\in\mathbb{N}:p|n\rightarrow p\leq X\},
\end{displaymath}
then $\alpha_{k}(n)$ is multiplicative, and $\alpha_{k}(n)=0$ if $n\notin S(X)$. We also have $\alpha_{k}(n)=d_{k}(n)$  if $n\in S(\sqrt{X})$, and $\alpha_{k}(p)=d_{k}(p)=k$. As in [\textbf{\ref{BK}}, Section 3], $|a_{k}(n)|\leq d_{3|k|/2}(n)$ and we can truncate the series, for $s=1/2$, at $D^{\theta}$,
\begin{equation}\label{8}
P_{X}^{*}(\chi_{8d})^{k}=\sum_{\substack{n\in S(X)\\n\leq D^\theta}}\frac{\alpha_{k}(n)\chi_{8d}(n)}{\sqrt{n}}+O(D^{-\delta\theta/4+\epsilon})
\end{equation}
Thus
\begin{equation}\label{9}
\sum_{0<d\leq D}{\!\!\!\!}^{\flat}\ P_{X}^{*}(\chi_{8d})^{k}=\sum_{\substack{n\in S(X)\\n\leq D^\theta}}\frac{\alpha_{k}(n)}{\sqrt{n}}\sum_{0<d\leq D}{\!\!\!\!}^{\flat}\ \chi_{8d}(n)+O(D^{1-\delta\theta/4+\epsilon}).
\end{equation}

We first consider the contribution of the main terms $n=m^2$ in the sum. By Lemma 2,
\begin{displaymath}
I=\frac{2D}{3\zeta(2)}\sum_{\substack{m\in S(X)\\2\nmid m\\m\leq D^{\theta/2}}}\frac{\alpha_{k}(m^2)a(m)}{m}+O\bigg(D^{1/2}\sum_{m\leq D^{\theta/2}}m^{-1+\epsilon}\bigg).
\end{displaymath}
The $O$-term is $\ll D^{1/2+\epsilon}$. As in \eqref{8}, the sum in the leading term can be extended to all $m\in S(X),2\nmid m$ with the gain of at most $O(D^{1-\theta/4+\epsilon})$. Since $\alpha_{k}(n)=d_{k}(n)$ for $n\in S(\sqrt{X})$, and $\alpha_{k}(p^2)=k^2/2$ for $\sqrt{X}<p\leq X$, the leading term of $I$ is
\begin{eqnarray*}
&&\frac{2D}{3\zeta(2)}\prod_{3\leq p\leq X}\bigg(\sum_{j=0}^{\infty}\frac{\alpha_{k}(p^{2j})a(p^j)}{p^j}\bigg)\\
&\sim&D^{*}\prod_{3\leq p\leq\sqrt{X}}\bigg(1+\sum_{j=1}^{\infty}\frac{d_{k}(p^{2j})}{p^{j-1}(p+1)}\bigg)\prod_{\sqrt{X}<p\leq X}\bigg(1+\frac{k^2}{2(p+1)}+O\bigg(\frac{1}{p^2}\bigg)\bigg)\\
&\sim&D^{*}\prod_{3\leq p\leq\sqrt{X}}\bigg(\bigg(1-\frac{1}{p}\bigg)^{k(k+1)/2}\bigg(1+\sum_{j=1}^{\infty}\frac{d_{k}(p^{2j})}{p^{j-1}(p+1)}\bigg)\bigg)\\
&&\qquad\qquad\qquad\qquad\prod_{3\leq p\leq\sqrt{X}}\bigg(1-\frac{1}{p}\bigg)^{-k(k+1)/2}\prod_{\sqrt{X}<p\leq X}\bigg(1-\frac{1}{p}\bigg)^{-k^2/2}\\
&\sim&2^{-k(k+2)/2}D^{*}(e^{\gamma}\log X)^{k(k+1)/2}\\
&&\qquad\qquad\prod_{p\leq\sqrt{X}}\bigg(\bigg(1-\frac{1}{p}\bigg)^{k(k+1)/2}\bigg(1+\sum_{j=1}^{\infty}\frac{d_{k}(p^{2j})}{p^{j-1}(p+1)}\bigg)\bigg).
\end{eqnarray*}
The product can be extended to include all primes $p$ as
\begin{displaymath}
\prod_{p>\sqrt{X}}\bigg(\bigg(1-\frac{1}{p}\bigg)^{k(k+1)/2}\bigg(1+\sum_{j=1}^{\infty}\frac{d_{k}(p^{2j})}{p^{j-1}(p+1)}\bigg)\bigg)=\prod_{p>\sqrt{X}}\bigg(1+O\bigg(\frac{1}{p^2}\bigg)\bigg)=1+o(1).
\end{displaymath}
So
\begin{displaymath}
I=(1+o(1))D^{*}a_{k,Sp}(e^{\gamma}\log X)^{k(k+1)/2}+O(D^{1-\theta/4+\epsilon}+D^{1/2+\epsilon}).
\end{displaymath}

The remaining terms in \eqref{9}, by the Polya-Vinogradov inequality, contribute to
\begin{displaymath}
J\ll\sum_{n\leq D^\theta}n^\epsilon\ll D^{\theta+\epsilon}.
\end{displaymath}
Choosing any $\theta<1$, and combining with Lemma 1, we obtain the theorem.

\subsection{Proof of Theorem S3.}

We divide the terms $d\leq D$ into dyadic blocks. Consider the block $D_{1}<d\leq 2D_{1}$. From \eqref{2} and \eqref{8} we have
\begin{eqnarray}\label{10}
\sum_{D_{1}<d\leq 2D_{1}}{\!\!\!\!\!\!\!\!}^{\flat}\ L({\scriptstyle{\frac{1}{2}}},\chi_{8d})P_{X}^{*}(\chi_{8d})^{-1}&=&\sum_{\substack{n\in S(X)\\n\leq D_{1}^\theta}}\frac{\alpha_{-1}(n)}{\sqrt{n}}\sum_{D_{1}<d\leq 2D_{1}}{\!\!\!\!\!\!\!\!}^{\flat}\ L({\scriptstyle{\frac{1}{2}}},\chi_{8d})\chi_{8d}(n)\nonumber\\
&&\qquad\qquad\qquad+O(D_{1}^{1-\delta\theta/4+\epsilon}).
\end{eqnarray}
Since $\alpha_{-1}(n)$ is supported on cube-free integers, and $\chi_{8d}(n)=0$ when $n$ is even, we can write $n=uv^2$, where $u$, $v$ are odd, square-free and $(u,v)=1$. From [\textbf{\ref{Sound1}}] (Proposition 1.1 and Proposition 1.2), we obtain
\begin{eqnarray*}
\sum_{D_{1}<d\leq 2D_{1}}{\!\!\!\!\!\!\!\!}^{\flat}\ L({\scriptstyle{\frac{1}{2}}},\chi_{8d})\chi_{8d}(n)&=&\\
&&\!\!\!\!\!\!\!\!\!\!\!\!\!\!\!\!\!\!\!\!\!\!\!\!\!\!\!\!\!\!\!\!\!\!\!\!\!\!\!\!\!\!\!\!\!\!\!\!\!\!\!\!(1+o(1))\frac{\sqrt{2}a_{1,Sp}D_{1}\log D_{1}}{3\zeta(2)}\frac{B_{1}(n)}{\sqrt{u}}-\frac{2\sqrt{2}a_{1,Sp}D_{1}}{3\zeta(2)}\frac{B_{1}(n)\log u}{\sqrt{u}}+E_1(n),
\end{eqnarray*}
where
\begin{displaymath}
B_{1}(n)=\prod_{p|n}\frac{p}{p+1}\bigg(1-\frac{1}{p(p+1)}\bigg)^{-1},
\end{displaymath}
and
\begin{displaymath}
\bigg|\sum_{n\leq D_{1}^\theta}\frac{\alpha_{-1}(n)}{\sqrt{n}}E_1(n)\bigg|\ll D_{1}^{3/4+\theta/2+\epsilon}.
\end{displaymath}
The first term above will contribute to the main term in \eqref{10} while the second term only contributes to an admissible error. We first consider the main term. Writing $P=\prod_{2<p\leq X}p$, for $n=uv^2\in S(X)$ we have $v|P$ and $u|(P/v)$. Noting that $B_{1}(n)$ is multiplicative, the main term in \eqref{10} is
\begin{eqnarray*}
I_1&=&\frac{\sqrt{2}a_{1,Sp}D_{1}\log D_{1}}{3\zeta(2)}\sum_{\substack{uv\in S(X)\\uv^2\leq D_{1}^\theta\\u,v\textrm{ odd}\\(u,v)=1}}\frac{\alpha_{-1}(uv^2)B_{1}(uv)}{uv}\\
&=&\frac{\sqrt{2}a_{1,Sp}D_{1}\log D_{1}}{3\zeta(2)}\sum_{\substack{v|P\\v\leq D_{1}^{\theta/2}}}\frac{\alpha_{-1}(v^2)B_{1}(v)}{v}\sum_{\substack{u|(P/v)\\u\leq D_{1}^\theta/v^2}}\frac{\alpha_{-1}(u)B_{1}(u)}{u}.
\end{eqnarray*}
We can extend both sums to all of $v|P$, and $u|(P/v)$, respectively, with the gain of at most `little \textit{o}' of the main term. Define the following multiplicative functions
\begin{displaymath}
T_{1}(n):=\sum_{h|n}\frac{\alpha_{-1}(h)B_{1}(h)}{h},
\end{displaymath}
and
\begin{displaymath}
T_{2}(n):=\sum_{h|n}\frac{\alpha_{-1}(h^2)B_{1}(h)}{hT_{1}(h)}.
\end{displaymath}
Then
\begin{eqnarray*}
I_1&\sim&\frac{\sqrt{2}a_{1,Sp}D_{1}\log D_{1}}{3\zeta(2)}T_{1}(P)T_{2}(P)\\
&\sim&\frac{\sqrt{2}a_{1,Sp}D_{1}\log D_{1}}{3\zeta(2)}\prod_{2<p\leq X}\bigg(T_{1}(p)+\frac{\alpha_{-1}(p^2)B_{1}(p)}{p}\bigg)\\
&\sim&\frac{\sqrt{2}a_{1,Sp}D_{1}\log D_{1}}{3\zeta(2)}\prod_{2<p\leq X}\bigg(1+\frac{\alpha_{-1}(p)B_{1}(p)}{p}+\frac{\alpha_{-1}(p^2)B_{1}(p)}{p}\bigg).
\end{eqnarray*}
We note that $\alpha_{-1}(p)=-1$, $\alpha_{-1}(p^2)=0$ if $p\leq\sqrt{X}$, and $\alpha_{-1}(p^2)=1/2$ if $\sqrt{X}<p\leq X$. Some straightforward calculations then give, if $p\leq\sqrt{X}$,
\begin{displaymath}
1+\frac{\alpha_{-1}(p)B_{1}(p)}{p}+\frac{\alpha_{-1}(p^2)B_{1}(p)}{p}=1-\frac{p}{p(p+1)-1},
\end{displaymath}
and, if $\sqrt{X}<p\leq X$, 
\begin{displaymath}
1+\frac{\alpha_{-1}(p)B_{1}(p)}{p}+\frac{\alpha_{-1}(p^2)B_{1}(p)}{p}=1-\frac{1}{2p}+O\bigg(\frac{1}{p^2}\bigg).
\end{displaymath}
Thus
\begin{eqnarray*}
I_1&\sim&\frac{\sqrt{2}a_{1,Sp}D_{1}\log D_{1}}{3\zeta(2)}\prod_{2<p\leq\sqrt{X}}\bigg(1-\frac{p}{p(p+1)-1}\bigg)\prod_{\sqrt{X}<p\leq X}\bigg(1-\frac{1}{2p}+O\bigg(\frac{1}{p^2}\bigg)\bigg)\\
&\sim&\frac{D_{1}\log D_{1}}{6\zeta(2)}\prod_{2<p\leq\sqrt{X}}\bigg(1-\frac{1}{p}\bigg)\prod_{\sqrt{X}<p\leq X}\bigg(1-\frac{1}{p}\bigg)^{1/2}\\
&\sim&\frac{\sqrt{2}D_{1}}{3\zeta(2)}\frac{\log D_{1}}{e^{\gamma}\log X}.
\end{eqnarray*}

The error is, as $B_{1}(n)\leq1$,
\begin{equation}\label{12}
J_1\ll D_{1}\sum_{uv\in S(X)}\frac{d(uv^2)\log u}{uv}\ll D_{1}\sum_{v\in S(X)}\frac{d(v^2)}{v}\sum_{u\in S(X)}\frac{d(u)\log u}{u}.
\end{equation}
We define
\begin{displaymath}
f(\sigma)=\sum_{u\in S(X)}\frac{d(u)}{u^\sigma}=\prod_{p\leq X}\frac{1}{(1-p^{-\sigma})^2}.
\end{displaymath}
Then
\begin{displaymath}
f'(\sigma)=-2f(\sigma)\sum_{p\leq X}\frac{\log p}{p^\sigma-1}.
\end{displaymath}
Thus the sum over $u$ is $-f'(1)\ll f(1)\sum_{p\leq X}\log p/p\ll(\log X)^3$. We also have $\sum_{v\in S(X)}d(v^2)/v\ll\prod_{p\leq X}(1-1/p)^{-3}\ll(\log X)^3$. So $J_1\ll D_{1}(\log X)^6$. Hence, choosing any $\theta<1/2$,
\begin{displaymath}
\sum_{D_{1}<d\leq 2D_{1}}{\!\!\!\!\!\!\!\!}^{\flat}\ L({\scriptstyle{\frac{1}{2}}},\chi_{8d})P_{X}^{*}(\chi_{8d})^{-1}=(1+o(1))\frac{\sqrt{2}D_{1}}{3\zeta(2)}\frac{\log D_{1}}{e^{\gamma}\log X}+O(D_{1}(\log X)^6).
\end{displaymath}
Summing over all the dyadic blocks and combining with Lemma 1 we obtain Theorem S3.

\subsection{Proof of Theorem S4.}

Consider the dyadic block $D_{1}<d\leq 2D_{1}$. From \eqref{3} and \eqref{8}, we have
\begin{eqnarray}\label{11}
\sum_{D_{1}<d\leq 2D_{1}}{\!\!\!\!\!\!\!\!}^{\flat}\ L({\scriptstyle{\frac{1}{2}}},\chi_{8d})^2P_{X}^{*}(\chi_{8d})^{-2}&=&\sum_{\substack{n\in S(X)\\n\leq D_{1}^\theta}}\frac{\alpha_{-2}(n)}{\sqrt{n}}\sum_{D_{1}<d\leq 2D_{1}}{\!\!\!\!\!\!\!\!}^{\flat}\ L({\scriptstyle{\frac{1}{2}}},\chi_{8d})^{2}\chi_{8d}(n)\nonumber\\
&&\qquad\qquad\qquad+O(D_{1}^{1-\delta\theta/4+\epsilon}).
\end{eqnarray}
As in [\textbf{\ref{BK}}] (cf. Lemma 7), we can assume that $\alpha_{-2}(n)$ is supported on cube-free integers. Since $\chi_{8d}(n)=0$ when $n$ is even, we can write $n=uv^2$, where $u$, $v$ are odd, square-free and $(u,v)=1$. Using Proposition 1.1 and Proposition 1.3 in [\textbf{\ref{Sound1}}], we have
\begin{eqnarray*}
\sum_{D_{1}<d\leq 2D_{1}}{\!\!\!\!\!\!\!\!}^{\flat}\ L({\scriptstyle{\frac{1}{2}}},\chi_{8d})^2\chi_{8d}(n)&=&(1+o(1))\frac{a_{2,Sp}D_{1}(\log D_{1})^3}{18\zeta(2)}\frac{d(u)B_{2}(uv)\sqrt{u}}{\sigma(u)}\\
&&\qquad+O\bigg(D_{1}(\log D_{1})^2\frac{d(u)\log u}{\sqrt{u}}\bigg)+E_{2}(n),
\end{eqnarray*}
where
\begin{displaymath}
B_{2}(n)=\prod_{p|n}\bigg(1+\frac{1}{p}+\frac{1}{p^2}-\frac{4}{p(p+1)}\bigg)^{-1},
\end{displaymath}
and
\begin{displaymath}
\bigg|\sum_{n\leq D_{1}^\theta}\frac{\alpha_{-2}(n)}{\sqrt{n}}E_2(n)\bigg|\ll D_{1}^{3/4+3\theta/4+\epsilon}.
\end{displaymath}
So we can write the main term in \eqref{11} as, say, $I_{2}+J_{2}$. We have
\begin{eqnarray*}
J_{2}&\ll&D_{1}(\log D_{1})^2\sum_{uv\in S(X)}\frac{d_{3}(uv^2)d(u)\log u}{uv}\\
&\ll&D_{1}(\log D_{1})^2\sum_{v\in S(X)}\frac{d_{3}(v^2)}{v}\sum_{u\in S(X)}\frac{d_{3}(u)d(u)\log u}{u}.
\end{eqnarray*}
As in \eqref{12}, the sum over $u$ is $\ll(\log X)^7$, and the sum over $v$ is $\ll(\log X)^6$. Hence $J_{2}\ll D_{1}(\log D_{1})^{2}(\log X)^{13}$. So the main contribution to \eqref{11} is
\begin{eqnarray*}
I_{2}&\sim&\frac{a_{2,Sp}D_{1}(\log D_{1})^3}{18\zeta(2)}\sum_{\substack{uv\in S(X)\\uv^2\leq D_{1}^\theta\\u,v\textrm{ odd}\\(u,v)=1}}\frac{\alpha_{-2}(uv^2)d(u)B_{2}(uv)}{\sigma(u)v}\\
&\sim&\frac{a_{2,Sp}D_{1}(\log D_{1})^3}{18\zeta(2)}\sum_{\substack{v|P\\v\leq D_{1}^{\theta/2}}}\frac{\alpha_{-2}(v^2)B_{2}(v)}{v}\sum_{\substack{u|(P/v)\\u\leq D_{1}^\theta/v^2}}\frac{\alpha_{-2}(u)d(u)B_{2}(u)}{\sigma(u)}.
\end{eqnarray*}
We can extend both sums to all of $v|P$, and $u|(P/v)$, respectively, with the gain of at most 'little \textit{o}' of the main term. Define the following multiplicative functions
\begin{displaymath}
T_{3}(n):=\sum_{h|n}\frac{\alpha_{-2}(h)d(h)B_{2}(h)}{\sigma(h)},
\end{displaymath}
\begin{displaymath}
T_{4}(n):=\sum_{h|n}\frac{\alpha_{-2}(h^2)B_{2}(h)}{hT_{3}(h)}.
\end{displaymath}
Then
\begin{eqnarray*}
I_{2}&\sim&\frac{a_{2,Sp}D_{1}(\log D_{1})^3}{18\zeta(2)}T_{3}(P)T_{4}(P)\\
&\sim&\frac{a_{2,Sp}D_{1}(\log D_{1})^3}{18\zeta(2)}\prod_{2<p\leq X}\bigg(T_{3}(p)+\frac{\alpha_{-2}(p^2)B_{2}(p)}{p}\bigg)\\
&\sim&\frac{a_{2,Sp}D_{1}(\log D_{1})^3}{18\zeta(2)}\prod_{2<p\leq X}\bigg(1+\frac{\alpha_{-2}(p)d(p)B_{2}(p)}{\sigma(p)}+\frac{\alpha_{-2}(p^2)B_{2}(p)}{p}\bigg).
\end{eqnarray*}
We note that $a_{2,Sp}\sim2^{-4}\prod_{2<p\leq X}(1-1/p)B_{2}(p)^{-1}$. So
\begin{displaymath}
I_{2}\sim\frac{D_{1}(\log D_{1})^3}{288\zeta(2)}\prod_{2<p\leq X}\bigg(1-\frac{1}{p}\bigg)\bigg(B_{2}(p)^{-1}+\frac{\alpha_{-2}(p)d(p)}{\sigma(p)}+\frac{\alpha_{-2}(p^2)}{p}\bigg).
\end{displaymath}
We have $\alpha_{-2}(p)=-2$, and $\alpha_{-2}(p^2)=1$ if $p\leq\sqrt{X}$, $\alpha_{-2}(p^2)=2$ if $\sqrt{X}<p\leq X$. Some straightforward calculations then give, if $p\leq\sqrt{X}$,
\begin{displaymath}
B_{2}(p)^{-1}+\frac{\alpha_{-2}(p)d(p)}{\sigma(p)}+\frac{\alpha_{-2}(p^2)}{p}=\bigg(1-\frac{1}{p}\bigg)^2,
\end{displaymath}
and, if $\sqrt{X}<p\leq X$, 
\begin{displaymath}
B_{2}(p)^{-1}+\frac{\alpha_{-2}(p)d(p)}{\sigma(p)}+\frac{\alpha_{-2}(p^2)}{p}=1-\frac{1}{p}+O\bigg(\frac{1}{p^2}\bigg).
\end{displaymath}
Thus
\begin{eqnarray*}
I_{2}&\sim&\frac{D_{1}(\log D_{1})^3}{288\zeta(2)}\prod_{2<p\leq\sqrt{X}}\bigg(1-\frac{1}{p}\bigg)^{3}\prod_{\sqrt{X}<p\leq X}\bigg(1-\frac{1}{p}\bigg)^{2}\\
&\sim&\frac{D_{1}}{18\zeta(2)}\bigg(\frac{\log D_{1}}{e^{\gamma}\log X}\bigg)^3.
\end{eqnarray*}
Choosing any $\theta<1/3$ and summing over all the dyadic blocks, the result follows.

\section{Orthogonal family}

\subsection{Proof of Theorem O2.}

We require some lemmas. We begin with a standard lemma [\textbf{\ref{G}}].
\begin{lemm}
For $m,n\geq1$,
\begin{displaymath}
\lambda_{f}(m)\lambda_{f}(n)=\sum_{\substack{d|(m,n)\\(d,q)=1}}\lambda_{f}\bigg(\frac{mn}{d^2}\bigg).
\end{displaymath}
Also if $p$ is a prime and $p\ne q$ then
\begin{eqnarray*}
\lambda_{f}(p)^{2m}=\sum_{r=0}^{m}\bigg(\binom{2m}{m-r}-\binom{2m}{m-r-1}\bigg)\lambda_{f}(p^{2r}),\\
\lambda_{f}(p)^{2m+1}=\sum_{r=0}^{m}\bigg(\binom{2m+1}{m-r}-\binom{2m+1}{m-r-1}\bigg)\lambda_{f}(p^{2r+1}).
\end{eqnarray*}
\end{lemm}

The next lemma is a particular case of Petersson's trace formula.
\begin{lemm}
For $m,n\geq1$, we have
\begin{displaymath}
\sum_{f}{\!}^{h}\ \lambda_{f}(m)\lambda_{f}(n)=\delta_{m,n}-J(m,n),
\end{displaymath}
where $\delta_{m,n}$ is the Kronecker symbol and
\begin{displaymath}
J(m,n)=2\pi\sum_{c\geq1}\frac{S(m,n;cq)}{cq}J_{1}\bigg(\frac{4\pi\sqrt{mn}}{cq}\bigg),
\end{displaymath}
in which $J_{1}(x)$ is the Bessel function of order $1$, and $S(m,n;c)$ is the Kloosterman sum
\begin{displaymath}
S(m,n;c)=\sum_{a\ (\emph{mod}\
c)}{\!\!\!\!\!\!\!}^{\displaystyle{*}}\ e\bigg(\frac{ma+n\overline{a}}{c}\bigg).
\end{displaymath}
Moreover we have
\begin{displaymath}
J(m,n)\ll(m,n,q)^{1/2}(mn)^{1/2+\epsilon}q^{-3/2}.
\end{displaymath}
\end{lemm}
The above estimate follows easily from the bound $J_{1}(x)\ll x$ and Weil's bound on Kloosterman sums.

We mention a result of Jutila [\textbf{\ref{J2}}] (cf. Theorem 1.7), which is an extension of the Voronoi summation formula.
\begin{lemm}
Let $f:\mathbb{R}^{+}\rightarrow\mathbb{C}$ be a $C^\infty$ function which vanishes in the neighbourhood of $0$ and is rapidly decreasing at infinity. Then for $c\geq1$ and $(a,c)=1$,
\begin{eqnarray*}
c\sum_{m\geq1}d(m)e\bigg(\frac{am}{c}\bigg)f(m)&=&2\int_{0}^{\infty}(\log\frac{\sqrt{x}}{c}+\gamma)f(x)dx\\
&&\ -2\pi\sum_{m\geq1}d(m)e\bigg(\frac{-\overline{a}m}{c}\bigg)\int_{0}^{\infty}Y_{0}\bigg(\frac{4\pi\sqrt{mx}}{c}\bigg)f(x)dx\\
&&\ +4\sum_{m\geq1}d(m)e\bigg(\frac{\overline{a}m}{c}\bigg)\int_{0}^{\infty}K_{0}\bigg(\frac{4\pi\sqrt{mx}}{c}\bigg)f(x)dx.
\end{eqnarray*}
\end{lemm}

\begin{lemm}
Let
\begin{displaymath}
P_{X}(f,s)=\exp\bigg(\sum_{n\leq X}\frac{\Lambda_{f}(n)}{n^{s}\log n}\bigg),
\end{displaymath}
so $P_{X}(f)=P_{X}(f,1/2)$, and let $P_{X,k}^{*}(f)=P_{X,k}^{*}(f,1/2)$, where
\begin{displaymath}
P_{X,k}^{*}(f,s)=\prod_{p\leq\sqrt{X}}\bigg(1-\frac{\lambda_{f}(p)}{p^s}+\frac{1}{p^{2s}}\bigg)^{-k}\prod_{\sqrt{X}<p\leq X}\bigg(1-\frac{|k|\lambda_{f}(p)}{p^s}+\frac{k^{2}\lambda_{f}(p)^{2}}{2p^{2s}}\bigg)^{-\emph{sign}(k)}.
\end{displaymath}
Then for any $k\in\mathbb{R}$ we have
\begin{displaymath}
P_{X}(f,s)^{k}=P_{X,k}^{*}(f,s)(1+O_{k}(X^{-1/6+\epsilon})),
\end{displaymath}
uniformly for $\sigma\geq1/2$.
\end{lemm}
\begin{proof}
Let $N_{p}=[\log X/\log p]$, the integer part of $\log X/\log p$. Then
\begin{eqnarray*}
P_{X}(f,s)^{k}&=&\exp\bigg(k\sum_{p\leq X}\sum_{1\leq j\leq N_{p}}\frac{\alpha_{f}(p)^{j}+\alpha'_{f}(p)^{j}}{jp^{js}}\bigg)\\
&=&\exp\bigg(k\sum_{p\leq\sqrt{X}}\sum_{1\leq j\leq N_{p}}\frac{\alpha_{f}(p)^{j}+\alpha'_{f}(p)^{j}}{jp^{js}}+k\sum_{\sqrt{X}<p\leq X}\frac{\lambda_{f}(p)}{p^{s}}\bigg),
\end{eqnarray*}
since $N_{p}=1$ for $\sqrt{X}<p\leq X$. We have
\begin{eqnarray*}
\bigg(1-\frac{|k|\lambda_{f}(p)}{p^s}+\frac{k^{2}\lambda_{f}(p)^{2}}{2p^{2s}}\bigg)^{-\textrm{sign}(k)}&=&\exp\bigg(\textrm{sign}(k)\sum_{j\geq1}\frac{1}{j}\bigg(\frac{|k|\lambda_{f}(p)}{p^s}-\frac{k^{2}\lambda_{f}(p)^{2}}{2p^{2s}}\bigg)^j\bigg)\\
&=&\exp\bigg(\frac{k\lambda_{f}(p)}{p^s}+O_{k}\bigg(\frac{1}{p^{3\sigma}}\bigg)\bigg).
\end{eqnarray*}
So
\begin{displaymath}
P_{X,k}^{*}(f,s)=\exp\bigg(k\sum_{p\leq \sqrt{X}}\sum_{j\geq1}\frac{\alpha_{f}(p)^{j}+\alpha'_{f}(p)^{j}}{jp^{js}}+k\sum_{\sqrt{X}<p\leq X}\bigg(\frac{\lambda_{f}(p)}{p^s}+O_{k}\bigg(\frac{1}{p^{3\sigma}}\bigg)\bigg)\bigg).
\end{displaymath}
Therefore
\begin{displaymath}
P_{X}(f,s)^{k}P_{X,k}^{*}(f,s)^{-1}=\exp\bigg(-k\sum_{p\leq\sqrt{X}}\sum_{j>N_{p}}\frac{\alpha_{f}(p)^{j}+\alpha'_{f}(p)^{j}}{jp^{js}}+O_{k}\bigg(\sum_{\sqrt{X}<p\leq X}\frac{1}{p^{3\sigma}}\bigg)\bigg).
\end{displaymath}
The argument in the exponent is
\begin{eqnarray*}
&\ll_k&\sum_{p\leq\sqrt{X}}\frac{1}{p^{\sigma(N_{p}+1)}}+\sum_{\sqrt{X}<p\leq X}\frac{1}{p^{3\sigma}}\\
&\ll_k&\sum_{2\leq j\leq\frac{\log X}{\log2}}\sum_{X^{1/(j+1)}<p\leq X^{1/j}}\frac{1}{p^{(j+1)/2}}+\sum_{\sqrt{X}<p\leq X}\frac{1}{p^{3/2}}\\
&\ll_k&\sum_{2\leq j\leq\frac{\log X}{\log2}}\frac{1}{X^{1/6}}+\frac{1}{X^{1/4}}\ll_k X^{-1/6+\epsilon}.
\end{eqnarray*}
Hence $P_{X}(f,s)^{k}P_{X,k}^{*}(f,s)^{-1}=1+O_{k}(X^{-1/6+\epsilon})$ as required.
\end{proof}

We proceed with the proof of the theorem. For $k\geq0$, the above lemma gives $P_{X}(f)^k=(1+o(1))P_{X,k}^{*}(f)$, where
\begin{displaymath}
P_{X,k}^{*}(f)=\prod_{p\leq\sqrt{X}}\bigg(1-\frac{\lambda_{f}(p)}{p^{1/2}}+\frac{1}{p}\bigg)^{-k}\prod_{\sqrt{X}<p\leq X}\bigg(1-\frac{k\lambda_{f}(p)}{p^{1/2}}+\frac{k^{2}\lambda_{f}(p)^{2}}{2p}\bigg)^{-1}.
\end{displaymath}
The second product is, using $\lambda_{f}(p)^2=\lambda_{f}(p^2)+1$,
\begin{eqnarray*}
&&\prod_{\sqrt{X}<p\leq X}\bigg(1+\frac{k\lambda_{f}(p)}{p^{1/2}}-\frac{k^{2}\lambda_{f}(p)^{2}}{2p}+\frac{k^2\lambda_{f}(p)^2}{p}+O_{k}\bigg(\frac{1}{p^{3/2}}\bigg)\bigg)\\
&=&(1+o(1))\prod_{\sqrt{X}<p\leq X}\bigg(\bigg(1+\frac{k^2}{2p}\bigg)+\frac{k\lambda_{f}(p)}{p^{1/2}}+\frac{k^2\lambda_{f}(p^2)}{2p}\bigg).
\end{eqnarray*}
We have, from Lemma 3,
\begin{eqnarray*}
\bigg(1-\frac{\lambda_{f}(p)}{p^{1/2}(1+1/p)}\bigg)^{-k}&=&\sum_{j\geq0}\frac{d_{k}(p^j)}{p^{j/2}(1+1/p)^j}\lambda_{f}(p)^j\\
&&\!\!\!\!\!\!\!\!\!\!\!\!\!\!\!\!\!\!\!\!\!\!\!\!\!\!\!\!\!\!\!\!\!\!\!\!\!\!\!\!\!\!\!\!\!\!\!\!\!\!\!\!\!\!\!\!\!\!\!\!\!\!=\sum_{j\geq0}\frac{d_{k}(p^{2j})}{p^{j}(1+1/p)^{2j}}\sum_{r=0}^{j}\bigg(\binom{2j}{j-r}-\binom{2j}{j-r-1}\bigg)\lambda_{f}(p^{2r})\\
&&\!\!\!\!\!\!\!\!\!\!\!\!\!\!\!\!\!\!\!\!\!\!\!\!\!\!\!\!\!\!\!\!\!\!\!\!\!\!\!\!\!\!\!\!\!\!\!\!\!\!\!\!\!\!\!\!\!\!\!\!\!\!\qquad+\sum_{j\geq0}\frac{d_{k}(p^{2j+1})}{p^{j+1/2}(1+1/p)^{2j+1}}\sum_{r=0}^{j}\bigg(\binom{2j+1}{j-r}-\binom{2j+1}{j-r-1}\bigg)\lambda_{f}(p^{2r+1})\\
&&\!\!\!\!\!\!\!\!\!\!\!\!\!\!\!\!\!\!\!\!\!\!\!\!\!\!\!\!\!\!\!\!\!\!\!\!\!\!\!\!\!\!\!\!\!\!\!\!\!\!\!\!\!\!\!\!\!\!\!\!\!\!=\sum_{r\geq0}\beta(p,r)\frac{\lambda_{f}(p^{r})}{p^{r/2}},
\end{eqnarray*}
where
\begin{displaymath}
\beta(p,2r)=\sum_{j\geq r}\frac{d_{k}(p^{2j})}{p^{j-r}(1+1/p)^{2j}}\bigg(\binom{2j}{j-r}-\binom{2j}{j-r-1}\bigg),
\end{displaymath}
and
\begin{displaymath}
\beta(p,2r+1)=\sum_{j\geq r}\frac{d_{k}(p^{2j+1})}{p^{j-r}(1+1/p)^{2j+1}}\bigg(\binom{2j+1}{j-r}-\binom{2j+1}{j-r-1}\bigg).
\end{displaymath}
So
\begin{eqnarray*}
P_{X}(f)^{k}&=&(1+o(1))\prod_{p\leq\sqrt{X}}\bigg(1+\frac{1}{p}\bigg)^{-k}\prod_{p\leq\sqrt{X}}\bigg(\sum_{r\geq0}\beta(p,r)\frac{\lambda_{f}(p^{r})}{p^{r/2}}\bigg)\\
&&\qquad\qquad\quad\prod_{\sqrt{X}<p\leq X}\bigg(\bigg(1+\frac{k^2}{2p}\bigg)+\frac{k\lambda_{f}(p)}{p^{1/2}}+\frac{k^2\lambda_{f}(p^2)}{2p}\bigg).
\end{eqnarray*}
Writing the above as $QT$ where
\begin{displaymath}
Q=\prod_{p\leq\sqrt{X}}\bigg(1+\frac{1}{p}\bigg)^{-k}\prod_{p\leq\sqrt{X}}\beta(p,0)\prod_{\sqrt{X}<p\leq X}\bigg(1+\frac{k^2}{2p}\bigg),
\end{displaymath}
and
\begin{displaymath}
T=\sum_{n\geq1}\gamma(n)\frac{\lambda_{f}(n)}{\sqrt{n}}.
\end{displaymath}
We have that $\gamma(n)$ is multiplicative, and $\gamma(n)=0$ if $n\notin S(X)$. Moreover, when $p\leq\sqrt{X}$, 
\begin{displaymath}
\gamma(p^j)=\frac{\beta(p,j)}{\beta(p,0)},
\end{displaymath}
and when $\sqrt{X}<p\leq X$,
\begin{displaymath}
\gamma(p^j)=\left\{ \begin{array}{ll}
\frac{k}{1+k^{2}/2p} &\qquad \textrm{if $j=1$}\\
\frac{k^2/2}{1+k^2/2p} &\qquad \textrm{if $j=2$}\\
0 & \qquad\textrm{if $j\geq3$.} 
\end{array} \right.
\end{displaymath}
We have
\begin{eqnarray}\label{38}
Q&\sim&\prod_{p\leq\sqrt{X}}\bigg(1-\frac{1}{p^2}\bigg)^{-k}\prod_{p\leq\sqrt{X}}\bigg(1-\frac{1}{p}\bigg)^{-k(k-1)/2}\nonumber\\
&&\qquad\prod_{\sqrt{X}<p\leq X}\bigg(1-\frac{1}{p}\bigg)^{-k^2/2}\prod_{p\leq\sqrt{X}}\bigg(1-\frac{1}{p}\bigg)^{k(k+1)/2}\beta(p,0)\nonumber\\
&\sim&a_{k,O}(e^{\gamma}\log X)^{k(k-1)/2},
\end{eqnarray}
as the last product can be extended to include all $p$ since $(1-1/p)^{k(k+1)/2}\beta(p,0)=1+O(1/p^2)$.

It is standard to check that $\gamma(p^j)<2^k j(\frac{2}{1+1/p})^j$. So $\gamma(n)\ll 2^{kw(n)}n\ll n^{1+\epsilon}$, where $w(n)$ is the number of distinct prime factors of $n$. Also, as in \eqref{8}, we can truncate the series $T$ at $q^\theta$ and have
\begin{displaymath}
T=\sum_{\substack{n\in S(X)\\n\leq q^\theta}}\gamma(n)\frac{\lambda_{f}(n)}{\sqrt{n}}+O(q^{-\delta\theta/4+\epsilon}).
\end{displaymath}
Applying Lemma 4, we obtain
\begin{eqnarray}\label{37}
\sum_{f}{\!}^{h}\ T&=&1+O\bigg(\sum_{n\leq q^\theta}|\gamma(n)|n^{\epsilon}q^{-3/2}\bigg)+O(q^{-\delta\theta/4+\epsilon})\nonumber\\
&=&1+O(q^{-3/2+2\theta+\epsilon}+q^{-\delta\theta/4+\epsilon}).
\end{eqnarray}
Choosing any $\theta<3/4$, the result follows from \eqref{38} and \eqref{37}.

\subsection{Proof of Theorems O3-O4.}

For $f\in S_{2}^{*}(q)$, and $1/2<c<1$, we consider
\begin{displaymath}
A_{k}(f,U):=\frac{1}{2\pi i}\int_{(c)}\frac{U^{s}\Lambda(f,s+{\scriptstyle{\frac{1}{2}}})^k}{(\frac{\sqrt{q}}{2\pi})^{k/2}}\frac{ds}{s},
\end{displaymath}
where $k=1$ or $k=2$. Moving the line of integration to $\Re s=-c$, and applying Cauchy's theorem and the functional equation, we derive that $A_{k}(f,U)=L(f,1/2)^{k}-\epsilon_{f}^{k}A_{k}(f, 1/U)$. Also, expanding $\Lambda(f,s+1/2)^k$ in a Dirichlet series and integrating termwise we get
\begin{equation}\label{5}
A_{1}(f,U):=\sum_{m\geq1}\frac{\lambda_{f}(m)}{\sqrt{m}}W_{1}\bigg(\frac{2\pi m}{\sqrt{q}U}\bigg),
\end{equation}
\begin{equation}\label{29}
A_{2}(f,U):=\sum_{m\geq1}\frac{d(m)\lambda_{f}(m)}{\sqrt{m}}W_{2}\bigg(\frac{4\pi^{2}m}{qU}\bigg),
\end{equation}
where
\begin{displaymath}
W_{1}(x)=\frac{1}{2\pi i}\int_{(c)}\Gamma(s+1)x^{-s}\frac{ds}{s},
\end{displaymath}
and
\begin{displaymath}
W_{2}(x)=\frac{1}{2\pi i}\int_{(c)}\Gamma(s+1)^{2}\zeta_{q}(2s+1)x^{-s}\frac{ds}{s}.
\end{displaymath}
We have $W_{k}(x)=O_{c}(x^{-c})$ and also, by moving the line of integration to $\Re s=-1+\epsilon$, $W_{1}(x)=1+O(x^{1-\epsilon})$, and
\begin{displaymath}
W_{2}(x)=-\bigg(1-\frac{1}{q}\bigg)\frac{\log x}{2}+\frac{\log q}{q}+O(x^{1-\epsilon}).
\end{displaymath}

For $k=1$, taking $U=q^\eta$, where $\eta$ is fixed and to be chosen later, we obtain $L(f,1/2)=A_{1}(f,q^\eta)+\epsilon_{f}A_{1}(f,q^{-\eta})$. We have $P_{X}(f)^{-1}=P_{X,-1}^{*}(f)(1+o(1))$. So
\begin{eqnarray}\label{17}
P_{X}(f)^{-1}&=&(1+o(1))\prod_{p\leq\sqrt{X}}\bigg(1-\frac{\lambda_{f}(p)}{\sqrt{p}}+\frac{1}{p}\bigg)\prod_{\sqrt{X}<p\leq X}\bigg(1-\frac{\lambda_{f}(p)}{\sqrt{p}}+\frac{\lambda_{f}(p)^{2}}{2p}\bigg)\nonumber\\
&=&(1+o(1))Q_{1}\prod_{p\leq\sqrt{X}}\bigg(1-\frac{1}{1+1/p}\frac{\lambda_{f}(p)}{\sqrt{p}}\bigg)\nonumber\\
&&\qquad\prod_{\sqrt{X}<p\leq X}\bigg(1-\frac{1}{1+1/2p}\frac{\lambda_{f}(p)}{\sqrt{p}}+\frac{1}{2(1+1/2p)}\frac{\lambda_{f}(p^2)}{p}\bigg),
\end{eqnarray}
where
\begin{equation}\label{23}
Q_{1}=\prod_{p\leq\sqrt{X}}\bigg(1+\frac{1}{p}\bigg)\prod_{\sqrt{X}<p\leq X}\bigg(1+\frac{1}{2p}\bigg)\sim\frac{e^{\gamma}\log X}{\sqrt{2}\zeta(2)}.
\end{equation}
Let us write the last two products in \eqref{17} as
\begin{displaymath}
T_{1}=\sum_{n\geq1}\frac{\beta_{1}(n)\lambda_{f}(n)}{\sqrt{n}}.
\end{displaymath}
We have that $\beta_{1}(n)$ is multiplicative, and $\beta_{1}(n)=0$ if $n\notin S(X)$. Moreover, when $p\leq\sqrt{X}$, 
\begin{displaymath}
\beta_{1}(p^j)=\left\{ \begin{array}{ll}
-\frac{1}{1+1/p} &\qquad \textrm{if $j=1$}\\
0 & \qquad\textrm{if $j\geq2$,} 
\end{array} \right.
\end{displaymath}
and when $\sqrt{X}<p\leq X$,
\begin{displaymath}
\beta_{1}(p^j)=\left\{ \begin{array}{ll}
-\frac{1}{1+1/2p} &\qquad \textrm{if $j=1$}\\
\frac{1}{2(1+1/2p)} &\qquad \textrm{if $j=2$}\\
0 & \qquad\textrm{if $j\geq3$.} 
\end{array} \right.
\end{displaymath}
It is easy to see that $|\beta_{1}(n)|\leq1$, so as in \eqref{8}, we can truncate the series at $q^\theta$ and obtain
\begin{equation}\label{20}
T_{1}=\sum_{\substack{n\in S(X)\\n\leq q^\theta}}\frac{\beta_{1}(n)\lambda_{f}(n)}{\sqrt{n}}+O(q^{-\delta\theta/4+\epsilon}).
\end{equation}

We write the truncated series as $R_1$. Using the bound on $W_{1}(x)$ and Lemma 4 we have
\begin{eqnarray}\label{18}
\sum_{f}{\!}^{h}\ \epsilon_{f}A_{1}(f,q^{-\eta})R_{1}&=&\sum_{\substack{m\geq1\\n\in S(X)\\n\leq q^\theta}}\frac{\beta_{1}(n)}{\sqrt{mn}}W_{1}\bigg(\frac{2\pi m}{q^{1/2-\eta}}\bigg)\sum_{f}{\!}^{h}\ \epsilon_{f}\lambda_{f}(m)\lambda_{f}(n)\nonumber\\
&\ll&q^{-1/2+\epsilon}\sum_{\substack{m\geq1\\n\leq q^\theta}}(mn)^{\epsilon}\bigg|W_{1}\bigg(\frac{2\pi m}{q^{1/2-\eta}}\bigg)\bigg|\ll q^{-\eta+\theta+\epsilon}.
\end{eqnarray}
Also
\begin{eqnarray}\label{19}
&&\sum_{f}{\!}^{h}\ A_{1}(f,q^\eta)R_{1}=\sum_{\substack{m\geq1\\n\in S(X)\\n\leq q^\theta}}\frac{\beta_{1}(n)}{\sqrt{mn}}W_{1}\bigg(\frac{2\pi m}{q^{1/2+\eta}}\bigg)\sum_{f}{\!}^{h}\ \lambda_{f}(m)\lambda_{f}(n)\nonumber\\
&=&\sum_{\substack{n\in S(X)\\n\leq q^\theta}}\frac{\beta_{1}(n)}{n}W_{1}\bigg(\frac{2\pi n}{q^{1/2+\eta}}\bigg)+O\bigg(q^{-3/2}\sum_{\substack{m\geq1\\n\leq q^\theta}}(mn)^{\epsilon}\bigg|W_{1}\bigg(\frac{2\pi m}{q^{1/2+\eta}}\bigg)\bigg|\bigg). 
\end{eqnarray}
The error term is $\ll q^{-1+\theta+\eta+\epsilon}$. The main term is
\begin{displaymath}
I_{1}=\sum_{\substack{n\in S(X)\\n\leq q^\theta}}\frac{\beta_{1}(n)}{n}+O\bigg(q^{-1/4-\eta/2}\sum_{n\leq q^\theta}\frac{1}{\sqrt{n}}\bigg)=\sum_{\substack{n\in S(X)\\n\leq q^\theta}}\frac{\beta_{1}(n)}{n}+O(q^{-1/4-\eta/2+\theta/2+\epsilon}).
\end{displaymath}
The sum, as in \eqref{8}, can be extended to all $n\in S(X)$ with the gain of at most $O(q^{-\theta/2+\epsilon})$. Thus, choosing any $0<\theta<\eta$ such that $\theta+\eta<1$,
\begin{eqnarray}\label{22}
I_{1}&=&\sum_{n\in S(X)}\frac{\beta_{1}(n)}{n}+O(q^{-\theta/2+\epsilon}+q^{-1/4-\eta/2+\theta/2+\epsilon})\nonumber\\
&\sim&\prod_{p\leq\sqrt{X}}\bigg(1+\frac{\beta_{1}(p)}{p}\bigg)\prod_{\sqrt{X}<p\leq X}\bigg(1+\frac{\beta_{1}(p)}{p}+\frac{\beta_{1}(p^2)}{p^2}\bigg)\nonumber\\
&\sim&\prod_{p\leq\sqrt{X}}\bigg(1-\frac{1}{p+1}\bigg)\prod_{\sqrt{X}<p\leq X}\bigg(1-\frac{1}{p+1/2}+O\bigg(\frac{1}{p^2}\bigg)\bigg)\nonumber\\
&\sim&\zeta(2)\prod_{p\leq\sqrt{X}}\bigg(1-\frac{1}{p}\bigg)\prod_{\sqrt{X}<p\leq X}\bigg(1-\frac{1}{p}\bigg)\sim\frac{\zeta(2)}{e^{\gamma}\log X}.
\end{eqnarray}
The result \eqref{13} and \eqref{17}-\eqref{22} together complete the proof for Theorem O3.\\

For $k=2$, we take $U=1$ and obtain $L(f,1/2)^{2}=2A_{2}(f,1)$. We have $P_{X}(f)^{-2}=P_{X,-2}^{*}(f)(1+o(1))$. So
\begin{eqnarray*}
P_{X}(f)^{-2}&=&(1+o(1))\prod_{p\leq\sqrt{X}}\bigg(1-\frac{\lambda_{f}(p)}{\sqrt{p}}+\frac{1}{p}\bigg)^{2}\prod_{\sqrt{X}<p\leq X}\bigg(1-\frac{2\lambda_{f}(p)}{\sqrt{p}}+\frac{2\lambda_{f}(p)^{2}}{p}\bigg)\\
&=&(1+o(1))\prod_{p\leq\sqrt{X}}\bigg(\bigg(1+\frac{3}{p}+\frac{1}{p^2}\bigg)-2\bigg(1+\frac{1}{p}\bigg)\frac{\lambda_{f}(p)}{\sqrt{p}}+\frac{\lambda_{f}(p^2)}{p}\bigg)\\
&&\qquad\prod_{\sqrt{X}<p\leq X}\bigg(\bigg(1+\frac{2}{p}\bigg)-\frac{2\lambda_{f}(p)}{\sqrt{p}}+\frac{2\lambda_{f}(p^2)}{p}\bigg),
\end{eqnarray*}
using $\lambda_{f}(p)^2=\lambda_{f}(p^2)+1$. Let us write the last two products as $Q_{2}T_{2}$ where 
\begin{displaymath}
Q_{2}=\prod_{p\leq\sqrt{X}}\bigg(1+\frac{3}{p}+\frac{1}{p^2}\bigg)\prod_{\sqrt{X}<p\leq X}\bigg(1+\frac{2}{p}\bigg),
\end{displaymath}
and
\begin{displaymath}
T_{2}=\sum_{n\geq1}\frac{\beta_{2}(n)\lambda_{f}(n)}{\sqrt{n}}.
\end{displaymath}
We note that $Q_{2}\ll(\log X)^3$. We have that $\beta_{2}(n)$ is multiplicative, and $\beta_{2}(n)=0$ if $n\notin S(X)$. Moreover, when $p\leq\sqrt{X}$, 
\begin{displaymath}
\beta_{2}(p^j)=\left\{ \begin{array}{ll}
-\frac{2(1+1/p)}{1+3/p+1/p^2} &\qquad \textrm{if $j=1$}\\
\frac{1}{1+3/p+1/p^2} &\qquad \textrm{if $j=2$}\\
0 & \qquad\textrm{if $j\geq3$,} 
\end{array} \right.
\end{displaymath}
and when $\sqrt{X}<p\leq X$,
\begin{displaymath}
\beta_{2}(p^j)=\left\{ \begin{array}{ll}
-\frac{2}{1+2/p} &\qquad \textrm{if $j=1$}\\
\frac{2}{1+2/p} &\qquad \textrm{if $j=2$}\\
0 & \qquad\textrm{if $j\geq3$.} 
\end{array} \right.
\end{displaymath}
It is easy to see that $|\beta_{2}(n)|\leq d(n)$, so as in \eqref{8}, we can truncate the series $T_2$ at $q^\theta$ and obtain
\begin{displaymath}
T_{2}=\sum_{\substack{n\in S(X)\\n\leq q^\theta}}\frac{\beta_{2}(n)\lambda_{f}(n)}{\sqrt{n}}+O(q^{-\delta\theta/4+\epsilon}).
\end{displaymath}
We have
\begin{eqnarray*}
\sum_{f}{\!}^{h}\ L(f,{\scriptstyle{\frac{1}{2}}})^{2}T_{2}&=&2\sum_{\substack{m\geq1\\n\in S(X)\\n\leq q^\theta}}\frac{d(m)\beta_{2}(n)}{\sqrt{mn}}W_{2}\bigg(\frac{4\pi^{2}m}{q}\bigg)\sum_{f}{\!}^{h}\ \lambda_{f}(m)\lambda_{f}(n)\\
&&\qquad\qquad\qquad\qquad\qquad+O(q^{-\delta\theta/4+\epsilon}). 
\end{eqnarray*}
We write this as $I_2+J_2$, say, where $I_2$ and $J_2$ are the diagonal and the off-diagonal contributions respectively.\\

We first consider $I_{2}$. We have
\begin{eqnarray*}
I_{2}&=&2\sum_{\substack{n\in S(X)\\n\leq q^\theta}}\frac{d(n)\beta_{2}(n)}{n}W_{2}\bigg(\frac{4\pi^2 n}{q}\bigg)\\
&=&(1+o(1))\sum_{\substack{n\in S(X)\\n\leq q^\theta}}\frac{d(n)\beta_{2}(n)}{n}\log\bigg(\frac{q}{n}\bigg)+O\bigg(\sum_{n\in S(X)}\frac{d(n)^2}{n}\bigg).
\end{eqnarray*}
The $O$-term is
\begin{displaymath}
\prod_{p\leq X}\bigg(1+\frac{d(p)^2}{p}+O\bigg(\frac{1}{p^2}\bigg)\bigg)\ll\prod_{p\leq X}\bigg(1-\frac{1}{p}\bigg)^{-4}\ll(\log X)^4.
\end{displaymath}
The error involving $\log n$ can be treated as in \eqref{12}, and is $\ll(\log X)^5$. So
\begin{displaymath}
I_{2}=(1+o(1))\log q\sum_{\substack{n\in S(X)\\n\leq q^{\theta}}}\frac{d(n)\beta_{2}(n)}{n}+O((\log X)^5).
\end{displaymath}
The sum, as in \eqref{8}, can be extended to all $n\in S(X)$ with the gain of at most $O(q^{-\theta/2+\epsilon})$. Hence the above expression is
\begin{eqnarray*}
&&(1+o(1))\log q\sum_{n\in S(X)}\frac{d(n)\beta_{2}(n)}{n}+O((\log X)^5)\\
&=&(1+o(1))\log q\prod_{p\leq X}\bigg(1+\frac{d(p)\beta_{2}(p)}{p}+\frac{d(p^2)\beta_{2}(p^2)}{p^2}\bigg)+O((\log X)^5)\\
&=&(1+o(1))\log q\prod_{p\leq\sqrt{X}}\bigg(1-\frac{4(1+1/p)}{p(1+3/p+1/p^2)}+\frac{3}{p^{2}(1+3/p+1/p^2)}\bigg)\\
&&\quad\qquad\prod_{\sqrt{X}<p\leq X}\bigg(1-\frac{4}{p(1+2/p)}+O\bigg(\frac{1}{p^2}\bigg)\bigg)+O((\log X)^5)\\
&=&(1+o(1))Q_{2}^{-1}\log q\prod_{p\leq\sqrt{X}}\bigg(1-\frac{1}{p}\bigg)\prod_{\sqrt{X}<p\leq X}\bigg(1-\frac{2}{p}\bigg)+O((\log X)^5)\\
&\sim&Q_{2}^{-1}\frac{\log q}{2e^{\gamma}\log X}.
\end{eqnarray*}

For $J_{2}$, using Lemma 4, we have
\begin{eqnarray*}
J_{2}&\ll&\sum_{\substack{m\geq1\\n\leq q^{\theta}}}d(m)d(n)\bigg|W_{2}\bigg(\frac{4\pi^2 m}{q}\bigg)\bigg|(mn)^{\epsilon}q^{-3/2}+q^{-\delta\theta/4+\epsilon}\\
&\ll&q^{-1/2+\theta+\epsilon}+q^{-\delta\theta/4+\epsilon}.
\end{eqnarray*}

The proof of Theorem O4 is completed by choosing any $\theta<1/2$.

\subsection{Proof of Theorem O5.}

Again, we first need a truncated series for $P_{X}(f)^{-3}$. From Lemma 6, we have
\begin{displaymath}
P_{X}(f)^{-3}=(1+o(1))\prod_{p\leq\sqrt{X}}\bigg(1-\frac{\lambda_{f}(p)}{\sqrt{p}}+\frac{1}{p}\bigg)^{3}\prod_{\sqrt{X}<p\leq X}\bigg(1-\frac{3\lambda_{f}(p)}{\sqrt{p}}+\frac{9\lambda_{f}(p)^{2}}{2p}\bigg).
\end{displaymath}
Using Lemma 3, the expression in the first product is
\begin{displaymath}
\bigg(1+\frac{6}{p}+\frac{6}{p^2}+\frac{1}{p^3}\bigg)-\bigg(3+\frac{8}{p}+\frac{3}{p^2}\bigg)\frac{\lambda_{f}(p)}{p^{1/2}}+\bigg(3+\frac{3}{p}\bigg)\frac{\lambda_{f}(p^2)}{p}-\frac{\lambda_{f}(p^3)}{p^{3/2}},
\end{displaymath}
and that in the second product is
\begin{displaymath}
\bigg(1+\frac{9}{2p}\bigg)-\frac{3\lambda_{f}(p)}{p^{1/2}}+\frac{9\lambda_{f}(p^2)}{2p}.
\end{displaymath}
So we can write $P_{X}(f)^{-3}$ as $Q_{3}T_{3}$, where
\begin{displaymath}
Q_{3}=\prod_{p\leq\sqrt{X}}\bigg(1+\frac{6}{p}+\frac{6}{p^2}+\frac{1}{p^3}\bigg)\prod_{\sqrt{X}<p\leq X}\bigg(1+\frac{9}{2p}\bigg),
\end{displaymath}
and
\begin{displaymath}
T_{3}=\sum_{n\geq1}\frac{\beta_{3}(n)\lambda_{f}(n)}{\sqrt{n}},
\end{displaymath}
where $\beta_{3}(n)$ is multiplicative and $\beta_{3}(n)=0$ if $n\notin S(X)$. Moreover, when $p\leq\sqrt{X}$, 
\begin{equation}\label{30}
\beta_{3}(p^j)=\left\{ \begin{array}{ll}
-\frac{3+8/p+3/p^2}{1+6/p+6/p^2+1/p^3} &\qquad \textrm{if $j=1$}\\
\frac{3+3/p}{1+6/p+6/p^2+1/p^3} &\qquad \textrm{if $j=2$}\\
-\frac{1}{1+6/p+6/p^2+1/p^3} &\qquad \textrm{if $j=3$}\\
0 & \qquad\textrm{if $j\geq4$,} 
\end{array} \right.
\end{equation}
and when $\sqrt{X}<p\leq X$,
\begin{equation}\label{31}
\beta_{3}(p^j)=\left\{ \begin{array}{ll}
-\frac{3}{1+9/2p} &\qquad \textrm{if $j=1$}\\
\frac{9}{2(1+9/2p)} &\qquad \textrm{if $j=2$}\\
0 & \qquad\textrm{if $j\geq3$.} 
\end{array} \right.
\end{equation}
It is easy to see that $|\beta_{3}(n)|\leq d_{3}(n)$, so as in \eqref{8}, we can truncate the series $T_3$ at $q^\theta$ and obtain
\begin{displaymath}
T_{3}=\sum_{\substack{n\in S(X)\\n\leq q^\theta}}\frac{\beta_{3}(n)\lambda_{f}(n)}{\sqrt{n}}+O(q^{-\delta\theta/4+\epsilon}).
\end{displaymath}

We write the truncated series as $R_3$. Letting $U=1$ in \eqref{29} and $U=q^{-\eta}$ in \eqref{5}, where $\eta$ is fixed and to be chosen later, we obtain $L(f,{\scriptstyle{\frac{1}{2}}})^3=M_{1}+\epsilon_{f}M_{2}$, where
\begin{equation}
M_1=2\sum_{l,m}\frac{d(l)}{\sqrt{lm}}W_{1}\bigg(\frac{2\pi m}{q^{1/2-\eta}}\bigg)W_{2}\bigg(\frac{4\pi^2l}{q}\bigg)\lambda_{f}(l)\lambda_{f}(m),
\end{equation}
and
\begin{equation}
M_2=2\sum_{l,m}\frac{d(l)}{\sqrt{lm}}W_{1}\bigg(\frac{2\pi m}{q^{1/2+\eta}}\bigg)W_{2}\bigg(\frac{4\pi^2l}{q}\bigg)\lambda_{f}(l)\lambda_{f}(m).
\end{equation}

\subsubsection{Contribution from $M_1$.}

From Lemma 3 and Lemma 4, we deduce that
\begin{equation}\label{49}
\sum_{f}{\!}^{h}\ M_{1}R_{3}=2\sum_{\substack{m\geq1\\uv\in S(X)\\uv\leq q^\theta}}\frac{\beta_{3}(uv)}{u\sqrt{mv}}W_{1}\bigg(\frac{2\pi mu}{q^{1/2-\eta}}\bigg)X_{1}(mv),
\end{equation}
where 
\begin{displaymath}
X_{1}(mv)=\sum_{l\geq1}\frac{d(l)}{\sqrt{l}}W_{2}\bigg(\frac{4\pi^2l}{q}\bigg)(\delta_{l,mv}-J(l,mv)).
\end{displaymath}
As $\theta$ and $\eta$ will be chosen to be small later ($3\theta+\eta<1/2$), throughout the proof we can restrict the sum over $m$ to $mv<q^{\Delta}$, for some $\Delta<1$. The contribution of the other terms, using the bound on $W_{1}$, is simply $\ll q^{-A}$ for every $A>0$. The bound on $J(l,mv)$ from Lemma 4 and the bound on $W_{2}(x)$ give
\begin{eqnarray*}
X_{1}(mv)&=&\frac{d(mv)}{\sqrt{mv}}W_{2}\bigg(\frac{4\pi^2mv}{q}\bigg)+O((mv)^{1/2+\epsilon}q^{-1/2+\epsilon})\\
&=&\frac{d(mv)}{2\sqrt{mv}}\log\frac{q}{4\pi^2mv}+O((mv)^{1/2+\epsilon}q^{-1/2+\epsilon}).
\end{eqnarray*}
The contribution of the $O$-term to \eqref{49} is bounded by
\begin{displaymath}
q^{-1/2+\epsilon}\sum_{\substack{m\leq q^{1/2-\eta}\\v\leq q^\theta}}(mv)^\epsilon\sum_{u\in S(X)}\frac{d_{3}(u)}{u}\ll q^{-\eta+\theta+\epsilon}.
\end{displaymath}
Thus
\begin{equation}\label{50}
\sum_{f}{\!}^{h}\ M_{1}R_{3}=\sum_{\substack{m\geq1\\uv\in S(X)\\uv\leq q^\theta}}\frac{\beta_{3}(uv)d(mv)}{umv}W_{1}\bigg(\frac{2\pi mu}{q^{1/2-\eta}}\bigg)\log\frac{q}{4\pi^2mv}+O(q^{-\eta+\theta+\epsilon}).
\end{equation}

\subsubsection{Contribution from $\epsilon_{f}M_2$.}

Applying Lemma 3 and Lemma 4 we have
\begin{equation}\label{39}
\sum_{f}{\!}^{h}\ \epsilon_{f}M_{2}R_{3}=2\sum_{\substack{m\geq1\\uv\in S(X)\\uv\leq q^\theta}}\frac{\beta_{3}(uv)}{u\sqrt{mv}}W_{1}\bigg(\frac{2\pi mu}{q^{1/2+\eta}}\bigg)X_{2}(mv),
\end{equation}
where
\begin{eqnarray}
X_{2}(mv)&=&-q^{1/2}\sum_{l\geq1}\frac{d(l)}{\sqrt{l}}W_{2}\bigg(\frac{4\pi^2l}{q}\bigg)J(ql,mv)\nonumber\\
&=&-\frac{2\pi}{\sqrt{q}}\sum_{c\geq1}\frac{1}{c}\sum_{l\geq1}\frac{d(l)}{\sqrt{l}}S(ql,mv;cq)J_{1}\bigg(\frac{4\pi}{c}\sqrt{\frac{lmv}{q}}\bigg)W_{2}\bigg(\frac{4\pi^2l}{q}\bigg).\label{51}
\end{eqnarray}

The sum over $c$ for which $q|c$, using Weil's bound for Kloosterman sums and $J_{1}(x)\ll x$, is
\begin{equation}\label{40}
\ll (mv)^{1/2+\epsilon}q^{-2+\epsilon}\sum_{l\leq q}l^\epsilon\ll q^{-1/2+\epsilon}.
\end{equation}

For $(c,q)=1$, we have
\begin{displaymath}
S(ql,mv;cq)=S(l,mv\overline{q};c)S(mv,0;q)=-S(l,mv\overline{q};c),
\end{displaymath}
since $S(mv,0;q)$ is a Ramanujan sum with $q$ prime and $(mv,q)=1$. Thus we need to study
\begin{equation}\label{41}
\frac{2\pi}{\sqrt{q}}\sum_{(c,q)=1}\frac{1}{c}\sum_{l\geq1}\frac{d(l)}{\sqrt{l}}S(l,mv\overline{q};c)J_{1}\bigg(\frac{4\pi}{c}\sqrt{\frac{lmv}{q}}\bigg)W_{2}\bigg(\frac{4\pi^2l}{q}\bigg).
\end{equation}
Using Weil's bound and the trivial bound $J_{1}(x)\ll x$, the tail of the series for $c>q^2$ is $O((mv)^{1/2+\epsilon}q^{-1+\epsilon})=O(q^{-1/2+\epsilon})$. We are therefore led to consider the terms $c<q^2$.

We fix a $C^{\infty}$ function $\xi:\mathbb{R}^{+}\rightarrow[0,1]$, which satisfies $\xi(x)=0$ for $0\leq x\leq1/2$ and $\xi(x)=1$ for $x\geq1$. We denote by $X_c$ the weighted inner sum in \eqref{41},
\begin{eqnarray*}
X_c&=&\sum_{l\geq1}\frac{d(l)}{\sqrt{l}}S(l,mv\overline{q};c)J_{1}\bigg(\frac{4\pi}{c}\sqrt{\frac{lmv}{q}}\bigg)W_{2}\bigg(\frac{4\pi^2l}{q}\bigg)\xi(l)\\
&=&\sum_{a\ (\textrm{mod}\
c)}{\!\!\!\!\!\!\!}^{\displaystyle{*}}\ e\bigg(\frac{mv\overline{qa}}{c}\bigg)\sum_{l\geq1}d(l)e\bigg(\frac{la}{c}\bigg)J_{1}\bigg(\frac{4\pi}{c}\sqrt{\frac{lmv}{q}}\bigg)W_{2}\bigg(\frac{4\pi^2l}{q}\bigg)\frac{\xi(l)}{\sqrt{l}}.
\end{eqnarray*}
Applying Lemma 5 we have
\begin{eqnarray*}
X_{c}&=&\frac{2}{c}S(mv,0;c)\int_{0}^{\infty}(\log\frac{\sqrt{x}}{c}+\gamma)J_{1}\bigg(\frac{4\pi}{c}\sqrt{\frac{xmv}{q}}\bigg)W_{2}\bigg(\frac{4\pi^2x}{q}\bigg)\frac{\xi(x)dx}{\sqrt{x}}\\
&&\!\!\!\!\!\!\!\!\!\!\!\!\!\!\!\!-\frac{2\pi}{c}\sum_{l\geq1}d(l)S(ql-mv,0;c)\int_{0}^{\infty}Y_{0}\bigg(\frac{4\pi\sqrt{lx}}{c}\bigg)J_{1}\bigg(\frac{4\pi}{c}\sqrt{\frac{xmv}{q}}\bigg)W_{2}\bigg(\frac{4\pi^2x}{q}\bigg)\frac{\xi(x)dx}{\sqrt{x}}\\
&&\!\!\!\!\!\!\!\!\!\!\!\!\!\!\!\!+\frac{4}{c}\sum_{l\geq1}d(l)S(ql+mv,0;c)\int_{0}^{\infty}K_{0}\bigg(\frac{4\pi\sqrt{lx}}{c}\bigg)J_{1}\bigg(\frac{4\pi}{c}\sqrt{\frac{xmv}{q}}\bigg)W_{2}\bigg(\frac{4\pi^2x}{q}\bigg)\frac{\xi(x)dx}{\sqrt{x}}.
\end{eqnarray*}
Using [\textbf{\ref{KM}}] (cf. Lemma 6 and Lemma 8), the contribution of the second term to $X_{2}(mv)$ is
\begin{equation}\label{45}
\ll(mv)^{1/2}q^{-1+\epsilon}\ll q^{-1/2+\epsilon},
\end{equation}
and that of the third term to $X_{2}(mv)$ is
\begin{equation}\label{46}
\ll_{\Delta}(mv)^{1/2}q^{-1+\epsilon}+q^{-1/2+\epsilon}\ll_{\Delta} q^{-1/2+\epsilon}.
\end{equation}

To deal with the contribution of the first term to $X_{2}(mv)$, we first remove the weight $\xi(x)$, which we can do with an admissible error
\begin{eqnarray*}
&&\frac{1}{\sqrt{q}}\sum_{c\leq q^2}\frac{1}{c^2}|S(mv,0;c)|\int_{0}^{1}\bigg|(\log\frac{\sqrt{x}}{c}+\gamma)J_{1}\bigg(\frac{4\pi}{c}\sqrt{\frac{xmv}{q}}\bigg)W_{2}\bigg(\frac{4\pi^2x}{q}\bigg)\bigg|\frac{dx}{\sqrt{x}}\\
&\ll&q^{-1/2}(\log q)^2\ll q^{-1/2+\epsilon}.
\end{eqnarray*}
It remains to study
\begin{eqnarray*}
&&\frac{4\pi}{\sqrt{q}}\sum_{\substack{c\leq q^2\\(c,q)=1}}\frac{1}{c^2}S(mv,0;c)\int_{0}^{\infty}(\log\frac{\sqrt{x}}{c}+\gamma)J_{1}\bigg(\frac{4\pi}{c}\sqrt{\frac{xmv}{q}}\bigg)W_{2}\bigg(\frac{4\pi^2x}{q}\bigg)\frac{dx}{\sqrt{x}}\\
&=&2\sum_{\substack{c\leq q^2\\(c,q)=1}}\frac{1}{c}S(mv,0;c)\int_{0}^{\infty}(\log\frac{\sqrt{qx}}{2\pi}+\gamma)J_{1}(2\sqrt{xmv})W_{2}(c^2x)\frac{dx}{\sqrt{x}}.
\end{eqnarray*}
Using the integral formula for $W_2(x)$, this is
\begin{equation}\label{47}
\frac{1}{2\pi i}\int_{(c)}2Z_{mv}(2s+1)\Gamma(s+1)^2\zeta_{q}(2s+1)L(s)\frac{ds}{s},
\end{equation}
where
\begin{displaymath}
Z_{mv}(s)=\sum_{\substack{c\leq q^2\\(c,q)=1}}\frac{S(mv,0;c)}{c^s},
\end{displaymath}
and
\begin{displaymath}
L(s)=\int_{0}^{\infty}(\log\frac{\sqrt{qx}}{2\pi}+\gamma)J_{1}(2\sqrt{xmv})x^{-s-1/2}dx.
\end{displaymath}
Applying [\textbf{\ref{KM}}] (cf. Lemma 3 and Lemma 4) and using the functional equation of $\Gamma(s)$, the integrand in \eqref{47} is
\begin{eqnarray*}
F(s)&=&\pi\bigg(\sigma_{-2s}(mv)\zeta_{q}(2s+1)^{-1}+O(q^{-4c+\epsilon})\bigg)\frac{\zeta_{q}(2s+1)}{\sin\pi s}\\
&&\qquad\bigg((mv)^{s-1/2}\bigg(\log\frac{q}{4\pi^2mv}+2\gamma+\frac{\Gamma'}{\Gamma}(1\pm s)\bigg)\bigg).
\end{eqnarray*}
The contribution of the $O$-term to the integral is easily seen to be $O(q^{-3c-1/2+\epsilon})$. The remaining term is an odd function of $s$. Moreover, the function is holomorphic in the strip $|\Re s|<1$, except for a single pole at $s=0$, and decreases exponentially in vertical strips. We denote by $X'_{2}(mv)$ the contribution of this to $X_{2}(mv)$. Shifting the contour to $\Re s=-c$, we deduce that 
\begin{eqnarray}\label{48}
2X'_{2}(mv)&=&\textrm{Res}_{s=0}\bigg(\pi\frac{\sigma_{-2s}(mv)(mv)^{s-1/2}}{\sin\pi x}\bigg(\log\frac{q}{4\pi^2mv}+2\gamma+\frac{\Gamma'}{\Gamma}(1\pm s)\bigg)\bigg)\nonumber\\
&=&\frac{d(mv)}{\sqrt{mv}}\log\frac{q}{4\pi^2mv}.
\end{eqnarray}

From \eqref{40}, \eqref{45}, \eqref{46} and \eqref{48}, we obtain that
\begin{displaymath}
X_{2}(mv)=\frac{d(mv)}{2\sqrt{mv}}\log\frac{q}{4\pi^2mv}+O(q^{-1/2+\epsilon}).
\end{displaymath}
The $O$-term contributes to \eqref{39} an error of size
\begin{displaymath}
q^{-1/2+\epsilon}\sum_{\substack{m\leq q^{1/2+\eta}\\v\leq q^\theta}}(mv)^{-1/2+\epsilon}\sum_{u\in S(X)}\frac{d_{3}(u)}{u}\ll q^{-1/4+\eta/2+\theta/2+\epsilon}.
\end{displaymath}
Thus
\begin{equation}\label{52}
\sum_{f}{\!}^{h}\ \epsilon_{f}M_{2}R_{3}=\sum_{\substack{m\geq1\\uv\in S(X)\\uv\leq q^\theta}}\frac{\beta_{3}(uv)d(mv)}{umv}W_{1}\bigg(\frac{2\pi mu}{q^{1/2+\eta}}\bigg)\log\frac{q}{4\pi^2mv}+O(q^{-1/4+\eta/2+\theta/2+\epsilon}).
\end{equation}

\subsubsection{Total contribution.}

From \eqref{50} and \eqref{52}, the main term of the third moment of $Z_{X}(f)$ is
\begin{displaymath}
M=Q_{3}\sum_{\substack{n\in S(X)\\n\leq q^\theta}}\frac{\beta_{3}(n)}{n}\sum_{\substack{m\geq1\\uv=n}}\frac{d(mv)}{m}\bigg[W_{1}\bigg(\frac{2\pi mu}{q^{1/2-\eta}}\bigg)+W_{1}\bigg(\frac{2\pi mu}{q^{1/2+\eta}}\bigg)\bigg]\log\frac{q}{4\pi^2mv}.
\end{displaymath}
We denote by $M^-$ and $M^+$ the sums involving $-\eta$ and $+\eta$ respectively. Using the integral formula for $W_1(x)$, the sum over $m,u$ and $v$ in $M^-$ is
\begin{eqnarray*}
&&\frac{1}{2\pi i}\int_{(c)}\Gamma(s+1)\sum_{uv=n}\sum_{m}\frac{d(mv)}{m}\bigg(\frac{q^{1/2-\eta}}{2\pi mu}\bigg)^{s}\log\frac{q}{4\pi^2mv}\frac{ds}{s}\\
&=&\frac{1}{2\pi i}\int_{(c)}\Gamma(s+1)\sum_{uv=n}\bigg(\frac{2\pi v}{uq^{1/2+\eta}}\bigg)^{s}\frac{\partial}{\partial s}\bigg[\bigg(\frac{q}{4\pi^2v}\bigg)^{s}\sum_{m}\frac{d(mv)}{m^{s+1}}\bigg]\frac{ds}{s}.
\end{eqnarray*}
We have
\begin{displaymath}
\sum_{m}\frac{d(mv)}{m^{s}}=\zeta(s)^2d(v)B_{v}(s),
\end{displaymath}
where
\begin{displaymath}
B_{v}(s)=\prod_{p^\alpha||v}\bigg(1-\frac{\alpha}{(\alpha+1)p^s}\bigg).
\end{displaymath}
Hence the integrand is
\begin{displaymath}
\Gamma(s+1)\frac{\zeta(s+1)^2}{s}\sum_{uv=n}d(v)\bigg(\frac{q^{1/2-\eta}}{2\pi u}\bigg)^{s}B_{v}(s+1)\bigg[\log\frac{q}{4\pi^2v}+2\frac{\zeta'}{\zeta}(s+1)+\frac{B'_{v}}{B_{v}}(s+1)\bigg].
\end{displaymath}
We shift the line of integration to $\Re s=-1/2$. On this line, as $X\ll(\log q)^{2-\delta}$,
\begin{displaymath}
B_{v}(s+1)\ll\prod_{p^\alpha||v}\bigg(1+\frac{1}{p^{1/2}}\bigg)\ll\prod_{p\leq X}\bigg(1+\frac{1}{p^{1/2}}\bigg)\ll q^{\epsilon},
\end{displaymath}
and
\begin{displaymath}
\frac{B'_{v}}{B_{v}}(s+1)\ll\sum_{p^\alpha||v}\frac{\log p}{p^{1/2}}\ll q^{\epsilon}.
\end{displaymath}
So the integral along $\Re s=-1/2$ is $\ll n^{3/2}q^{-1/4+\eta/2+\epsilon}$. Thus the main contribution to the contour integral comes from the multiple poles at $s=0$, which, with respect to $q$, gives a polynomial in $\log q$. To find the leading term, we can replace all $\zeta(s+1)$ factors by $1/s$ and so obtain
\begin{eqnarray*}
&&\sum_{uv=n}d(v)B_{v}(1)\textrm{Res}_{s=0}\frac{1}{s^3}\bigg(\frac{q^{1/2-\eta}}{2\pi u}\bigg)^{s}\bigg(\log q-\frac{2}{s}+O(\log n)\bigg)\\
&=&\frac{(1/2-\eta)^2(1+\eta)}{3}\sum_{uv=n}d(v)B_{v}(1)((\log q)^3+O((\log q)^2\log n)).
\end{eqnarray*}

Similarly to $M^+$, adding up $M^-$ and $M^+$ we obtain
\begin{eqnarray}\label{53}
M&=&\frac{Q_3}{6}\sum_{\substack{n\in S(X)\\n\leq q^{\theta}}}\frac{\beta_{3}(n)}{n}\sum_{uv=n}d(v)B_{v}(1)((\log q)^3+O((\log q)^2\log n))\nonumber\\
&&\qquad\qquad\qquad+O(q^{-1/4+\eta/2+3\theta/2+\epsilon}).
\end{eqnarray}
Let $I_3$ be the main term in $M$. Again we can extend the sum over $n$ to all $n\in S(X)$, with the gain of at most $O(q^{-\theta/2+\epsilon})$. We note that the function $f(n)=\sum_{uv=n}d(v)B_{v}(1)$ is multiplicative. Thus
\begin{displaymath}
I_{3}=\frac{(\log q)^3}{6}Q_{3}\prod_{p\leq X}\bigg(\sum_{j\geq0}\frac{\beta_{3}(p^j)}{p^j}\sum_{uv=p^j}d(v)B_{v}(1)\bigg)+O(q^{-\theta/2+\epsilon}).
\end{displaymath}
The sum over $u,v$ is
\begin{equation}\label{34}
\frac{(j+1)(j+2)}{2}-\frac{j(j+1)}{2p}.
\end{equation}
Using \eqref{30} and \eqref{31}, after some calculations we have that the expression in the bracket is, for $p\leq\sqrt{X}$,
\begin{displaymath}
\bigg(1+\frac{6}{p}+\frac{6}{p^2}+\frac{1}{p^3}\bigg)^{-1}\bigg(1-\frac{1}{p}\bigg)^3,
\end{displaymath}
and, for $\sqrt{X}<p\leq X$, is
\begin{displaymath}
\bigg(1+\frac{9}{2p}\bigg)^{-1}\bigg(1-\frac{9}{2p}+O\bigg(\frac{1}{p^2}\bigg)\bigg).
\end{displaymath}
Thus
\begin{eqnarray}\label{32}
I_{3}&\sim&\frac{(\log q)^3}{6}\prod_{p\leq\sqrt{X}}\bigg(1-\frac{1}{p}\bigg)^{3}\prod_{\sqrt{X}<p\leq X}\bigg(1-\frac{9}{2p}\bigg)\nonumber\\
&\sim&\frac{1}{12\sqrt{2}}\bigg(\frac{\log q}{e^\gamma\log X}\bigg)^3.
\end{eqnarray}

We now need to show that
\begin{equation}\label{33}
J_{3}=Q_{3}\sum_{n\in S(X)}\frac{d_{3}(n)}{n}\sum_{uv=n}d(v)B_{v}(1)\log n=o\bigg(\frac{\log q}{(\log X)^3}\bigg).
\end{equation}
We note from \eqref{34} that $\sum_{uv=n}d(v)B_{v}(1)\ll d_{3}(n)$, so
\begin{displaymath}
J_{3}\ll (\log X)^{6}\sum_{n\in S(X)}\frac{d_{3}(n)^{2}\log n}{n}.
\end{displaymath}
As in \eqref{12}, the sum is $\ll (\log X)^{10}$, so $J_{3}\ll(\log X)^{16}$.

Choosing any $0<\theta<\eta$ such that $3\theta+\eta<1/2$, the theorem now follows from \eqref{32} and \eqref{33}.

\subsection{Proof of Theorem O6.}

From Lemma 6, we have
\begin{displaymath}
P_{X}(f)^{-4}=(1+o(1))\prod_{p\leq\sqrt{X}}\bigg(1-\frac{\lambda_{f}(p)}{\sqrt{p}}+\frac{1}{p}\bigg)^{4}\prod_{\sqrt{X}<p\leq X}\bigg(1-\frac{4\lambda_{f}(p)}{\sqrt{p}}+\frac{8\lambda_{f}(p)^{2}}{p}\bigg).
\end{displaymath}
Using Lemma 3, the expression in the first product is
\begin{eqnarray*}
&&\bigg(1+\frac{10}{p}+\frac{20}{p^2}+\frac{10}{p^3}+\frac{1}{p^4}\bigg)-\bigg(4+\frac{20}{p}+\frac{20}{p^2}+\frac{4}{p^3}\bigg)\frac{\lambda_{f}(p)}{p^{1/2}}\\
&&\qquad\qquad+\bigg(6+\frac{15}{p}+\frac{6}{p^2}\bigg)\frac{\lambda_{f}(p^2)}{p}-\bigg(4+\frac{4}{p}\bigg)\frac{\lambda_{f}(p^3)}{p^{3/2}}+\frac{\lambda_{f}(p^4)}{p^{2}},
\end{eqnarray*}
and that in the second product is
\begin{displaymath}
\bigg(1+\frac{8}{p}\bigg)-\frac{4\lambda_{f}(p)}{p^{1/2}}+\frac{8\lambda_{f}(p^2)}{p}.
\end{displaymath}
So we can write the above as $Q_{4}T_{4}$, where
\begin{displaymath}
Q_{4}=\prod_{p\leq\sqrt{X}}\bigg(1+\frac{10}{p}+\frac{20}{p^2}+\frac{10}{p^3}+\frac{1}{p^4}\bigg)\prod_{\sqrt{X}<p\leq X}\bigg(1+\frac{8}{p}\bigg),
\end{displaymath}
and
\begin{displaymath}
T_{4}=\sum_{n\geq1}\frac{\beta_{4}(n)\lambda_{f}(n)}{\sqrt{n}}.
\end{displaymath}
Here $\beta_{4}(n)$ is multiplicative, and $\beta_{4}(n)=0$ if $n\notin S(X)$. Moreover, when $p\leq\sqrt{X}$, 
\begin{equation}\label{24}
\beta_{4}(p^j)=\left\{ \begin{array}{ll}
-\frac{4+20/p+20/p^2+4/p^3}{1+10/p+20/p^2+10/p^3+1/p^4} &\qquad \textrm{if $j=1$}\\
\frac{6+15/p+6/p^2}{1+10/p+20/p^2+10/p^3+1/p^4} &\qquad \textrm{if $j=2$}\\
-\frac{4+4/p}{1+10/p+20/p^2+10/p^3+1/p^4} &\qquad \textrm{if $j=3$}\\
\frac{1}{1+10/p+20/p^2+10/p^3+1/p^4} &\qquad \textrm{if $j=4$}\\
0 & \qquad\textrm{if $j\geq5$,} 
\end{array} \right.
\end{equation}
and when $\sqrt{X}<p\leq X$,
\begin{equation}\label{25}
\beta_{4}(p^j)=\left\{ \begin{array}{ll}
-\frac{4}{1+8/p} &\qquad \textrm{if $j=1$}\\
\frac{8}{1+8/p} &\qquad \textrm{if $j=2$}\\
0 & \qquad\textrm{if $j\geq3$.} 
\end{array} \right.
\end{equation}
It is easy to see that $|\beta_{4}(n)|\leq d_{4}(n)$, so as in \eqref{8}, we can truncate the series $T_4$ at $q^\theta$ and obtain
\begin{displaymath}
T_{4}=\sum_{\substack{n\in S(X)\\n\leq q^\theta}}\frac{\beta_{4}(n)\lambda_{f}(n)}{\sqrt{n}}+O(q^{-\delta\theta/4+\epsilon}).
\end{displaymath}
For $n\leq q^\theta$, where $\theta<4/33$, Theorem 1.2 and Theorem 4.1 in [\textbf{\ref{KMV1}}] give
\begin{eqnarray}
\sum_{f}{\!}^{h}\ L(f,{\scriptstyle{\frac{1}{2}}})^{4}\lambda_{f}(n)&=&\frac{2}{\sqrt{n}}\sum_{uv=n}\sum_{m\geq1}\frac{d(mu)d(mv)}{m}W_{2}\bigg(\frac{4\pi^2mu}{q}\bigg)\log\frac{q}{4\pi^2mv}\nonumber\\
&&\qquad\qquad\qquad+M^{OOD}(n)+O(n^{3/4}q^{-1/12+\epsilon})\label{54}.
\end{eqnarray}
The off-off-diagonal term is given by
\begin{displaymath}
M^{OOD}(n)=\frac{4}{\sqrt{n}}\frac{1}{(2\pi i)^2}\int_{(1.7)}\int_{(0.6)}I(n;s,t)\frac{ds}{s}\frac{dt}{t},
\end{displaymath}
where
\begin{eqnarray*}
I(n;s,t)&=&\Gamma(1+s)\Gamma(1-s)\Gamma(1+t)\Gamma(1-t)\prod\zeta(1\pm s\pm t)\\
&&\!\!\!\!\!\!\!\!\!\!\!\!\!\!\!\!\!\!\!\!\!\!\!\!\!\!\!\!\!\times\frac{1}{n^t}\sum_{aAD|n}A^{2t}\frac{\varphi(A)}{A}D^{t}\frac{\sigma_{2s}(D)}{D^s}\frac{\varphi(aD)}{aD}\frac{\mu(a)}{a}\sum_{(w,aAD)=1}\frac{\mu(w)}{w^2}\\
&&\!\!\!\!\!\!\!\!\!\!\!\!\!\!\!\!\!\!\!\!\!\!\!\!\!\!\!\!\!\times\bigg[\bigg(\log q-\log\frac{nw^2}{D}+2\sum_{p|aD}\frac{\log p}{p-1}-2\gamma+\frac{\Gamma'}{\Gamma}(1\pm s)\bigg)\\
&&\!\!\!\!\!\!\!\!\!\!\!\!\!\!\!\!\!\!\!\!\!\!\!\!\!\!\!\!\!\quad\times\bigg(\log q-\log Dw^2+2\sum_{p|A}\frac{\log p}{p-1}-2\gamma+\frac{\Gamma'}{\Gamma}(1\pm t)\bigg)\\
&&\!\!\!\!\!\!\!\!\!\!\!\!\!\!\!\!\!\!\!\!\!\!\!\!\!\!\!\!\!\quad+\bigg(2\log q-\log nw^4+2\sum_{p|aD}\frac{\log p}{p-1}+2\sum_{p|A}\frac{\log p}{p-1}-4\gamma\\
&&\!\!\!\!\!\!\!\!\!\!\!\!\!\!\!\!\!\!\!\!\!\!\!\!\!\!\!\!\!\quad\quad\quad+\frac{\Gamma'}{\Gamma}(1\pm s)+\frac{\Gamma'}{\Gamma}(1\pm t)\bigg)\bigg(-2\log2\pi+\sum\frac{\zeta'}{\zeta}(1\pm s\pm t)\bigg)\\
&&\!\!\!\!\!\!\!\!\!\!\!\!\!\!\!\!\!\!\!\!\!\!\!\!\!\!\!\!\!\quad-\frac{2\pi^2\sin^{2}\pi t}{(\cos\pi t+\cos\pi s)^2}+4(\log2\pi)^3-4\log2\pi\sum\frac{\zeta'}{\zeta}(1\pm s\pm t)\\
&&\!\!\!\!\!\!\!\!\!\!\!\!\!\!\!\!\!\!\!\!\!\!\!\!\!\!\!\!\!\quad+\sum\frac{\zeta''}{\zeta}(1\pm s\pm t)+2\sum{\!}^{'}\frac{\zeta'}{\zeta}(1\pm s\pm t)\frac{\zeta'}{\zeta}(1\pm s\pm t)\bigg].
\end{eqnarray*}
We note that $I(n;s,t)$ is even in both $s$ and $t$. Thus
\begin{displaymath}
M^{OOD}(n)=\frac{1}{4}\textrm{Res}_{s=t=0}\frac{4I(n;s,t)}{\sqrt{n}st}=\textrm{Res}_{t=0}\frac{I(n;0,t)}{\sqrt{n}t}.
\end{displaymath}
As we are only interested in the leading term, which is expected to involve two factors of $\log q$, we are allowed to simplify the expression for $I(n;0,t)$ by dropping terms like $\Gamma(1\pm t),\sum\log p/(p-1),\gamma,\Gamma'/\Gamma(1\pm t),\ldots$ and replacing $\zeta(1\pm t)$ by $\pm1/t$. We are thus led to study
\begin{eqnarray*}
&&\frac{1}{\sqrt{n}}\textrm{Res}_{t=0}\frac{1}{t^5}\sum_{aAD|n}\bigg(\frac{A^2D}{n}\bigg)^{t}\frac{\varphi(A)}{A}d(D)\frac{\varphi(aD)}{aD}\frac{\mu(a)}{a}\sum_{(w,aAD)=1}\frac{\mu(w)}{w^2}\\
&&\qquad\qquad\bigg[\bigg(\log q-\log\frac{nw^2}{D}\bigg)\bigg(\log q-\log Dw^2\bigg)+\frac{4}{t^2}\bigg],
\end{eqnarray*}
where the term $4/t^2$ comes from the sums of $\zeta''/\zeta$ and $(\zeta'/\zeta)(\zeta'/\zeta)$. This is bounded by
\begin{displaymath}
\frac{(\log q)^2(\log n)^4}{\sqrt{n}}\sum_{aAD|n}\frac{d(D)}{a}\ll\frac{(\log q)^2d(n)^4(\log n)^4}{\sqrt{n}}.
\end{displaymath}
Thus the contribution of $M^{OOD}(n)$ and the $O$-term in \eqref{54} to $\sum_{f}^{h}L(f,{\scriptstyle{\frac{1}{2}}})^{4}P_{X}(f)^{-4}$ is at most
\begin{displaymath}
Q_{4}(\log q)^2\sum_{n\in S(X)}\frac{d_{4}(n)d(n)^4(\log n)^4}{n}+O(q^{-1/12+3\theta/4+\epsilon}).
\end{displaymath}
As in \eqref{12}, this is
\begin{equation}\label{55}
\ll(\log X)^{78}(\log q)^2+O(q^{-1/12+3\theta/4+\epsilon}).
\end{equation}

So the leading term of $\sum_{f}^{h}L(f,{\scriptstyle{\frac{1}{2}}})^{4}P_{X}(f)^{-4}$ is expected to come from
\begin{equation}\label{36}
Q_{4}\sum_{\substack{n\in S(X)\\n\leq q^{\theta}}}\frac{2\beta_{4}(n)}{n}\sum_{uv=n}\sum_{m\geq1}\frac{d(mu)d(mv)}{m}W_{2}\bigg(\frac{4\pi^2mu}{q}\bigg)\log\frac{q}{4\pi^2mv}.
\end{equation}
Using the integral formula for $W_2(x)$, the sum over $m,u$ and $v$ is
\begin{eqnarray*}
&&\frac{1}{2\pi i}\int_{(c)}\Gamma(s+1)^{2}\zeta_{q}(2s+1)\sum_{uv=n}\sum_{m\geq1}\frac{d(mu)d(mv)}{m}\bigg(\frac{q}{4\pi^2mu}\bigg)^{s}\log\frac{q}{4\pi^2mv}\frac{ds}{s}\\
&=&\frac{1}{2\pi i}\int_{(c)}\Gamma(s+1)^{2}\zeta_{q}(2s+1)\sum_{uv=n}\bigg(\frac{v}{u}\bigg)^{s}\frac{\partial}{\partial s}\bigg[\bigg(\frac{q}{4\pi^2v}\bigg)^{s}\sum_{m\geq1}\frac{d(mu)d(mv)}{m^{s+1}}\bigg]\frac{ds}{s}.
\end{eqnarray*}
We have
\begin{displaymath}
\sum_{m\geq1}\frac{d(mu)d(mv)}{m^{s}}=\frac{\zeta^4(s)}{\zeta(2s)}d(u)d(v)B_{u,v}(s),
\end{displaymath}
where
\begin{displaymath}
B_{u,v}(s)=\prod_{\substack{p^\alpha||u\\p^\beta||v}}\bigg(1-\frac{(3\alpha\beta+2\alpha+2\beta)p^s-\alpha\beta}{(\alpha+1)(\beta+1)p^s(p^s+1)}\bigg).
\end{displaymath}
Hence the integrand is
\begin{eqnarray*}
&&\Gamma(s+1)^{2}\zeta_{q}(2s+1)\frac{\zeta^4(s+1)}{s\zeta(2(s+1))}\sum_{uv=n}d(u)d(v)\bigg(\frac{q}{4\pi^2u}\bigg)^{s}B_{u,v}(s+1)\\
&&\qquad\bigg[\log\frac{q}{4\pi^2v}+4\frac{\zeta'}{\zeta}(s+1)-2\frac{\zeta'}{\zeta}(2(s+1))+\frac{B'_{u,v}}{B_{u,v}}(s+1)\bigg].
\end{eqnarray*}
We shift the line of integration to $\Re s=-1/2+\epsilon$. It is easy to check that on this line, as $X\ll(\log q)^{2-\delta}$,
\begin{displaymath}
B_{u,v}(s+1)\ll\prod_{\substack{p^\alpha||u\\p^\beta||v}}\bigg(1+\frac{3}{p^{1/2+\epsilon}}\bigg)\ll\prod_{p\leq X}\bigg(1+\frac{3}{p^{1/2+\epsilon}}\bigg)\ll q^{\epsilon},
\end{displaymath}
and
\begin{displaymath}
\frac{B'_{u,v}}{B_{u,v}}(s+1)\ll\sum_{\substack{p^\alpha||u\\p^\beta||v}}\frac{\log p}{p^{1/2+\epsilon}}\ll q^{\epsilon}.
\end{displaymath}
So the integral along $\Re s=-1/2+\epsilon$ is $\ll n^{3/2}q^{-1/2+\epsilon}$. Thus the main contribution to the contour integral comes from the multiple poles at $s=0$, which, with respect to $q$, gives a polynomial in $\log q$. To find the leading term, we can replace all $\zeta(1+s)$ factors by $1/s$ and so obtain
\begin{eqnarray*}
&&\frac{1}{2\zeta(2)}\sum_{uv=n}d(u)d(v)B_{u,v}(1)\textrm{Res}_{s=0}\frac{1}{s^6}\bigg(\frac{q}{4\pi^2u}\bigg)^{s}\bigg(\log q-\frac{4}{s}+O(\log n)\bigg)\\
&=&\frac{1}{6!\zeta(2)}\sum_{uv=n}d(u)d(v)B_{u,v}(1)((\log q)^6+O((\log q)^5\log n)).
\end{eqnarray*}

We first consider the contribution of the main term above to \eqref{36}, which is
\begin{displaymath}
I_4=\frac{(\log q)^6}{360\zeta(2)}Q_{4}\sum_{\substack{n\in S(X)\\n\leq q^{\theta}}}\frac{\beta_{4}(n)}{n}\sum_{uv=n}d(u)d(v)B_{u,v}(1).
\end{displaymath}
Again we extend the sum over $n$ to all $n\in S(X)$, with the gain of at most $O(q^{-\theta/2+\epsilon})$. We note that the function $f(n)=\sum_{uv=n}d(u)d(v)B_{u,v}(1)$ is multiplicative. Thus
\begin{displaymath}
I_4=\frac{(\log q)^6}{360\zeta(2)}Q_{4}\prod_{p\leq X}\bigg(\sum_{j\geq0}\frac{\beta_{4}(p^j)}{p^j}\sum_{uv=p^j}d(u)d(v)B_{u,v}(1)\bigg)+O(q^{-\theta/2+\epsilon}).
\end{displaymath}
The sum over $u,v$ is
\begin{equation}\label{28}
(j+1)+j(j+1)\frac{1-1/p}{1+1/p}+\frac{(j-1)j(j+1)}{6}\frac{(1-1/p)^2}{1+1/p}.
\end{equation}
Using \eqref{24} and \eqref{25}, after some calculations we have that the expression in the bracket, for $p\leq\sqrt{X}$, is
\begin{eqnarray*}
&&\bigg(1+\frac{10}{p}+\frac{20}{p^2}+\frac{10}{p^3}+\frac{1}{p^4}\bigg)^{-1}\bigg(1+\frac{1}{p}\bigg)^{-1}\bigg(1-\frac{5}{p}+\frac{10}{p^2}-\frac{10}{p^3}+\frac{5}{p^4}-\frac{1}{p^5}\bigg)\\
&=&\bigg(1+\frac{10}{p}+\frac{20}{p^2}+\frac{10}{p^3}+\frac{1}{p^4}\bigg)^{-1}\bigg(1-\frac{1}{p^2}\bigg)^{-1}\bigg(1-\frac{1}{p}\bigg)^{6},
\end{eqnarray*}
and, for $\sqrt{X}<p\leq X$, is
\begin{displaymath}
\bigg(1+\frac{8}{p}\bigg)^{-1}\bigg(1-\frac{8}{p}+O\bigg(\frac{1}{p^2}\bigg)\bigg).
\end{displaymath}
Thus
\begin{eqnarray}\label{26}
I_{4}&\sim&\frac{(\log q)^6}{360\zeta(2)}\prod_{p\leq\sqrt{X}}\bigg(1-\frac{1}{p^2}\bigg)^{-1}\bigg(1-\frac{1}{p}\bigg)^{6}\prod_{\sqrt{X}<p\leq X}\bigg(1-\frac{8}{p}\bigg)\nonumber\\
&\sim&2^{-8}\frac{(\log q)^6}{360}\bigg(\frac{2}{e^\gamma\log X}\bigg)^{6}\sim\frac{1}{1440}\bigg(\frac{\log q}{e^\gamma\log X}\bigg)^6.
\end{eqnarray}

We now need to show that
\begin{equation}\label{27}
J_{4}=Q_{4}\sum_{n\in S(X)}\frac{d_{4}(n)}{n}\sum_{uv=n}d(u)d(v)B_{u,v}(1)\log n=o\bigg(\frac{\log q}{(\log X)^6}\bigg).
\end{equation}
We note from \eqref{28} that $\sum_{uv=n}d(u)d(v)B_{u,v}(1)\ll d_{4}(n)$, so
\begin{displaymath}
J_{4}\ll (\log X)^{10}\sum_{n\in S(X)}\frac{d_{4}(n)^{2}\log n}{n}.
\end{displaymath}
As in \eqref{12}, the sum is $\ll (\log X)^{17}$, so $J_{4}\ll(\log X)^{27}$.

Choosing any $\theta<1/9$, the theorem now follows from \eqref{55}, \eqref{26} and \eqref{27}.

\specialsection*{\textbf{Acknowledgements}}
We are grateful to the referee for pointing out an error in our original proof of Theorem O6, and for several helpful suggestions. JPK is supported by an EPSRC Senior Research Fellowship.

\end{document}